\documentclass[10pt]{article}
\usepackage{amsfonts}
\usepackage{hyperref}
\usepackage{mathrsfs,enumerate,amssymb,amsmath,color,lineno}
\usepackage{indentfirst}
\usepackage{amsmath}
\usepackage{amssymb}
\usepackage[amsmath,thmmarks]{ntheorem}
\usepackage{cite}
\usepackage{graphicx}
\usepackage{tikz}
\usepackage{subfigure}

\newtheorem{definition}{\bf Definition}[section]
\newtheorem{lemma}{\bf Lemma}[section]
\newtheorem{theorem}{\bf Theorem}[section]
\newtheorem{remark}{\bf Remark}[section]
\newtheorem{corollary}{\bf Corollary}[section]
\newtheorem{example}{\bf Example}[section]

\newtheorem{proposition}{\bf Proposition}[section]

\begin{document}
\setcounter{page}{1}

\title{{\textbf{Associativity of two-place functions generated by left continuous monotone functions and other properties}}\thanks {Supported by
the National Natural Science Foundation of China (No.12471440)}}
\author{Meng Chen\footnote{\emph{E-mail address}: mathchen2019@163.com}, Xue-ping Wang\footnote{Corresponding author. xpwang1@hotmail.com; fax: +86-28-84761502},\\
\emph{School of Mathematical Sciences, Sichuan Normal University,}\\
\emph{Chengdu 610066, Sichuan, People's Republic of China}}

\newcommand{\pp}[2]{\frac{\partial #1}{\partial #2}}
\date{}
\maketitle
\begin{quote}
{\bf Abstract} This article introduces a weak pseudo-inverse of a monotone function, which is applied to characterize the associativity of a two-place function $T: [0,1]^2\rightarrow [0,1]$ defined by $T(x,y)=t^{[-1]}(F(t(x),t(y)))$ where $F:[0,\infty]^2\rightarrow[0,\infty]$ is an associative function with neutral element in $[0,\infty]$, $t: [0,1]\rightarrow [0,\infty]$ is a left continuous monotone function and $t^{[-1]}:[0,\infty]\rightarrow[0,1]$ is the weak pseudo-inverse of $t$. It shows that the associativity of the function $T$ depends only on properties of the range of $t$. Moreover, it investigates the idempotence, the limit property, the conditional cancellation law and the continuity of the function $T$, respectively.

{\textbf{\emph{Keywords}}:} Monotone function; Weak pseudo-inverse; Left continuous function; Associative function; Triangular norm\\
\end{quote}

\section{Introduction}
\begin{figure}[htbp]
	\centering
\subfigure[] {
\includegraphics[width=6cm]{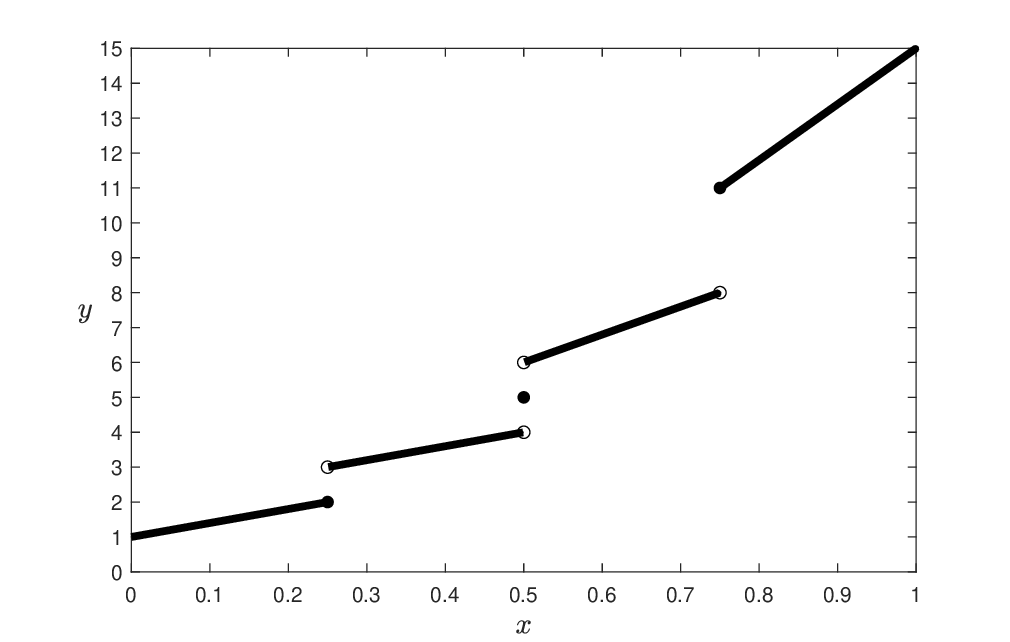}
	}
	\quad
	\label{T1}
\subfigure[]{
		\includegraphics[width=6cm]{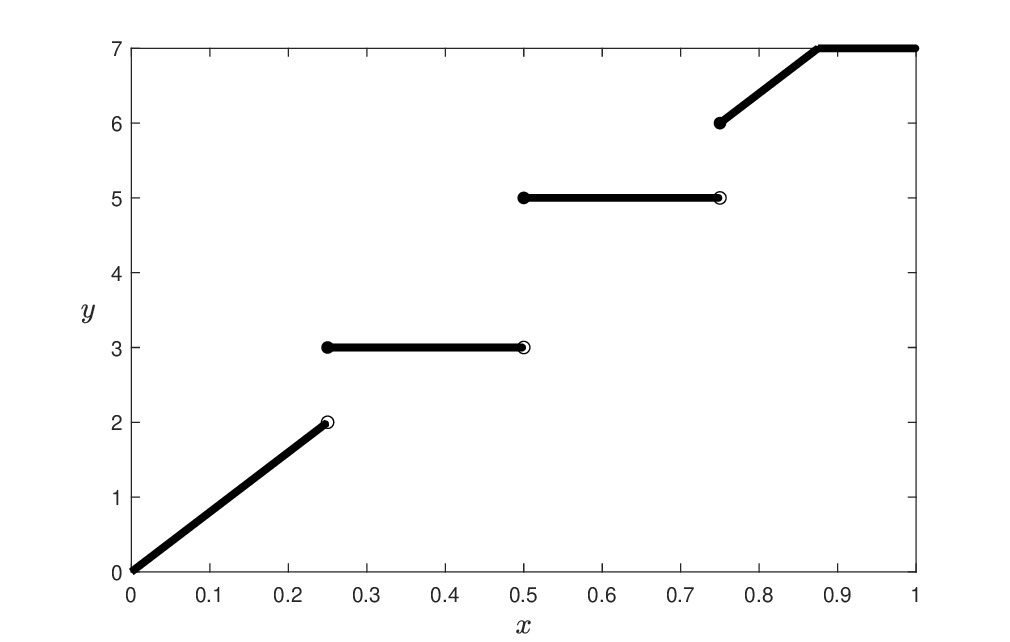}
	}
	\label{T2}
\end{figure}

In 1826, Abel \cite{Abel} obtained an easily checking result: Let $t:\mbox{Dom}(t)\rightarrow R$ with $\mbox{Dom}(t)\subseteq R$ be a continuous strictly monotone function whose range is closed under addition. Then the two-place function $T:(\mbox{Dom}(t))^{2}\rightarrow \mbox{Dom}(t)$ defined by
\begin{equation}
\label{eq:1}
T(x,y)=g(t(x)+t(y)),
\end{equation}
where $\mbox{Dom}(t)$ is the domain of $t$, $g:\mbox{Dom}(t)\rightarrow \mbox{Dom}(t)$ is the inverse function of $t$ and $R$ is the set of all real numbers, is associative. This result can be seen as the starting point of constructing a two-place real function that has nice algebraic properties through a monotone one-place real function (Note, the idea goes back to Abel \cite{Abel}). Following this idea, Schweizer and Sklar \cite{BS1961} and Ling \cite{CH1965} constructed triangular norms (t-norms for short) by continuous strictly decreasing functions, respectively. In particular, Klement, Mesiar and Pap \cite{EP2000} defined an additive generator of a t-norm $T$ as a strictly decreasing function $t :[0,1]\rightarrow [0,\infty]$ that is right continuous at $0$ with $t(1) = 0$ such that for all $(x,y)\in [0,1]^2$,
\begin{equation}
\label{eq:1}
t(x)+t(y)\in \mbox{Ran}(t)\cup [t(0),\infty],
\end{equation}
and
\begin{equation}
\label{eq:2}
 T(x,y)= g(t(x)+t(y))
\end{equation}
where $g$ is a pseudo-inverse of $t$ and $\mbox{Ran}(t)$ is a range of $t$, and they further pointed out that we can generalize the additive generator of a t-norm $T$ as it just satisfies \eqref{eq:2}. This idea was identified by Vicen\'{\i}k \cite{PV2005} when $t$ is a strictly monotone function shown in Figure \ref{T1}(a) where
$t$ is strictly increasing. The related work can refer to \cite{PV1998,PV1998b,PV2008} also. Zhang and Wang \cite{YM2024} proved that a right continuous monotone function, shown in Figure \ref{T2} (b) in which only the right continuous non-decreasing function is presented, may be a generator of an associative two-place function (see Corollaries 5.2 and 5.3 of \cite{YM2024}). Recently, Chen, Zhang and Wang \cite{2025chen} showed that a large number of monotone functions, shown in Figure \ref{T1} (c) in which only the non-decreasing function is presented, can be a generator of an associative two-place function. It is clear that they greatly generalized the related work of \cite{YO2008}, \cite{PV2005} and \cite{YM2024}. One naturally wishes that a left continuous monotone function, shown in Figure \ref{T1} (d) in which only the left continuous non-decreasing function is presented, can be a generator of an associative two-place function. However, Remark 4.1 of Zhang and Wang \cite{YM2024} showed that generally, \eqref{eq:2} is not associative when $t$ is a left continuous monotone function and $g$ is a pseudo-inverse of $t$. In this article, we consider a problem: what is the characterization of a left continuous monotone function that generates an associative two-place function?

The rest of this article is organized as follows. In Section 2, we mainly recall some basic concepts. In Section 3, we define a weak pseudo-inverse of a monotone function, and develop its properties. In Section 4, we give a representation of the range $\mbox{Ran}(t)$ of a left continuous non-decreasing function $t$. In Section 5, we first define an operation $\otimes$ on the $\mbox{Ran}(t)$, investigate some necessary and sufficient conditions for the operation $\otimes$ being associative and characterize what properties of $\mbox{Ran}(t)$ are equivalent to the associativity of a two-place function $T$ generated by a left continuous non-decreasing function $t$. Section 6 is devoted to exploring the idempotence, the limit property, the conditional cancellation law and the continuity of the function $T$, respectively. A conclusion is drawn in Section 7.

\begin{figure}[htbp]
	\centering
	\subfigure[]{
		\includegraphics[width=6cm]{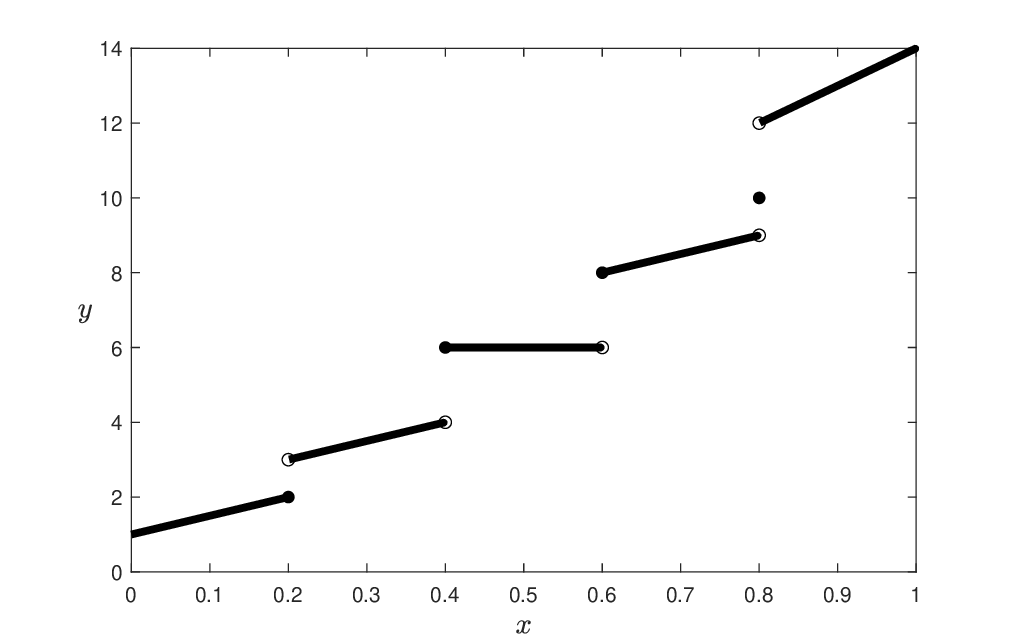}
	}
	\quad
	\label{T2}
\subfigure[]{
		\includegraphics[width=6cm]{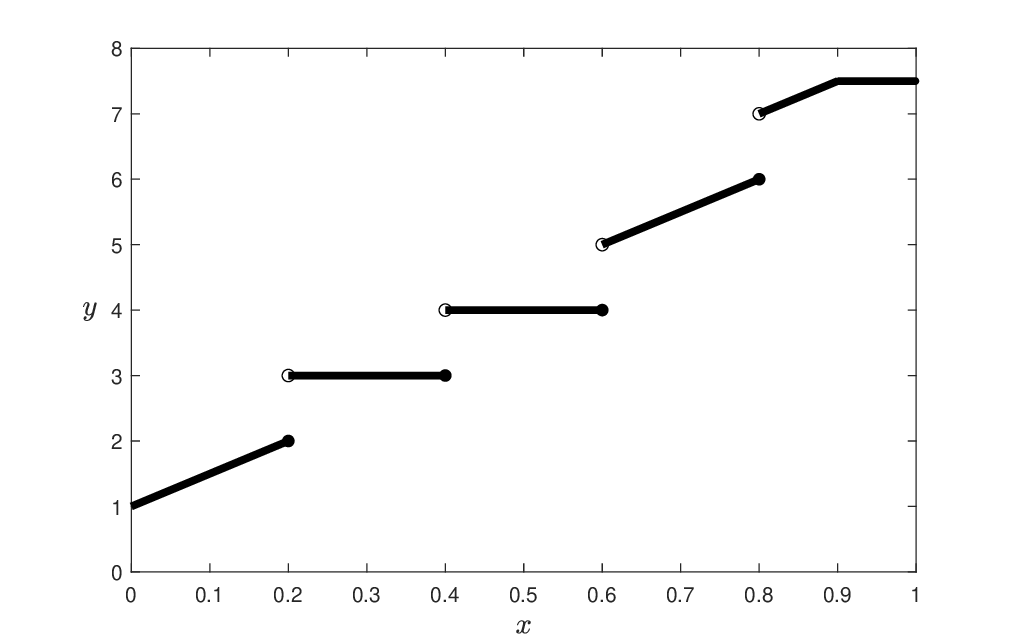}
	}
\caption{(a) a strictly increasing function; (b) a right continuous non-decreasing function; (c) a kind of non-decreasing functions; (d) a left continuous non-decreasing function.}
\end{figure}

\section{Preliminaries}

In this section, we recall some known basic concepts and results that will be used latter.
\begin{definition}[\cite{EP2000}]\label{de2.2}
\emph{A t-norm is a binary operator $T:[0, 1]^2\rightarrow [0, 1]$ such that for all $x, y, z\in[0, 1]$ the following conditions are satisfied:}

$(T1)$  $T(x,y)=T(y,x)$,

$(T2)$  $T(T(x,y),z)=T(x,T(y,z))$,

\emph{$(T3)$  $T(x,y)\leq T(x,z)$ whenever $y\leq z$,}

$(T4)$  $T(x,1)=x$.
\end{definition}

A binary operator $T:[0, 1]^2\rightarrow [0, 1]$ is called a t-subnorm if it satisfies $(T1), (T2), (T3)$, and $T(x,y)\leq \min\{x,y\}$ for all $x,y\in [0, 1]$.

\begin{definition}[\cite{EP2000}]\label{de2.2}
\emph{A t-conorm is a binary operator $S:[0, 1]^2\rightarrow [0, 1]$ such that for all $x, y, z\in[0, 1]$ the following conditions are satisfied:}

$(S1)$ $S(x,y)=S(y,x)$,

$(S2)$  $S(S(x,y),z)=S(x,S(y,z))$,

\emph{$(S3)$  $S(x,y)\leq S(x,z)$ whenever $y\leq z$,}

$(S4)$ $S(x,0)=x$.
\end{definition}

A binary operator $S:[0, 1]^2\rightarrow [0, 1]$ is called a t-supconorm if it satisfies $(S1), (S2), (S3)$, and $S(x,y)\geq \max\{x,y\}$ for all $x,y\in [0, 1]$.

\begin{definition}[\cite{EP2000,PV2005}]\label{def2.3}
\emph{Let $a, b, m, n\in [-\infty, \infty]$ with $a<b, m<n$ and $t:[a,b]\rightarrow[m,n]$ be a monotone function. Then the function $t^{(-1)}:[m,n]\rightarrow[a,b]$ defined by
\begin{equation*}
t^{(-1)}(y)=\sup\{x\in [a,b]\mid (t(x)- y)(t(b)-t(a))<0 \}
\end{equation*}
is called a pseudo-inverse of the monotone function $t$.}
\end{definition}

Let $[a,b]\subseteq [-\infty,\infty]$ with $a\leq b$. Then by convention, $\sup \emptyset=a$ and $\inf \emptyset=b$.
\begin{definition}[\cite{EP1999}]
\emph{Let $a, b, m, n\in [-\infty, \infty]$ with $a<b, m<n$ and $t:[a,b]\rightarrow[m,n]$ be a monotone non-decreasing function. Then each function $t^* :[m,n] \rightarrow[a, b]$ satisfying
\renewcommand{\labelenumi}{(\roman{enumi})}
\begin{enumerate}
\item $t\circ t^* \circ t =t$,
\item $t^{\wedge}\leq t^*\leq t^{\vee}$,
\end{enumerate}
is called a quasi-inverse of $t$, where functions $t^{\wedge}:[m,n] \rightarrow[a, b]$ and $t^{\vee}:[m,n] \rightarrow[a, b]$ are defined
by, respectively, $t^{\wedge}(y)= \sup t^{-1}([m,y))$,  $t^{\vee}(y)= \inf t^{-1}((y,n])$ in which $t^{-1}$ is an inverse function of $t$.}
\end{definition}

\begin{theorem}[\cite{EP2000}]\label{th2.1}
 Let $t :[0,1]\rightarrow [0,\infty]$ be a strictly decreasing function with $t(1) = 0$ such that
\begin{equation*}
t(x)+t(y)\in \emph{Ran}(t)\cup [t(0^{+}),\infty]
\end{equation*}
for all $(x,y)\in[0,1]^{2}$. Then the function $T: [0,1]^2\rightarrow [0,1]$ given by
\begin{equation*}
 T(x,y)= t^{(-1)}(t(x)+ t(y))
\end{equation*}
is a t-norm.
\end{theorem}

We write $N=\{1,2,\cdots, n, \cdots\}$.
\begin{definition}[\cite{EP2000}]\label{def2.2}
\emph{A binary function $T\colon[0,1]^{2}\to[0,1]$ is continuous if for all (non-decreasing/non-increasing) sequences $\{x_n\}_{n\in N}, \{y_n\}_{n\in N}, x_n, y_n\in [0,1]$ with $\lim_ {n \rightarrow \infty} x_n=x$ and $\lim_ {n \rightarrow \infty} y_n=y$, we have \begin{equation*}
\lim_ {n \rightarrow \infty} T(x_n, y_n)=T(x, y).
\end{equation*}}
\end{definition}

\begin{definition}
Let a function $F:[0,\infty]^2\rightarrow[0,\infty]$ be such that $([0,\infty],F,\leq)$
is a fully ordered Abel semigroup with $F(x,0) \geq x$ for all $x \in [0,\infty]$. If $F(x,y)<F(x,z)$ whenever $x<\infty$ and $y<z$ for all $x,y,z\in[0,\infty]$ then
$F$ is called to be strictly monotone. If $F$ is continuous and strictly monotone then $F$ is called to be strict.
\end{definition}
\begin{proposition}\label{pr1.6}
	Let a function $F:[0,\infty]^2\rightarrow[0,\infty]$ be such that $([0,\infty],F,\leq)$
is a fully ordered Abel semigroup with $F(x,0) \geq x$ for all $x \in [0,\infty]$. If $F$ is continuous, Then $a\in[0,\infty]$ is an idempotent element of $F$ if and only if   $F(a,x)=\max\{a,x\}$ for all $x\in[0,\infty]$.
\end{proposition}
\begin{proof}If for all  $x\in[0,\infty]$, $F(a,x)=\max\{a,x\}$, then, obviously, $F(a,a)=a$. Conversely, If $F(a,a)=a$, then for all $x\in[0,a]$, we have $$a\leq F(a,0)\leq F(a,x)\leq F(a,a)=a,$$ i.e., $F(a,x)=a$.  On the other hand, Because of $F(a,a)=a$ and $F(a,\infty)=\infty$,  the  continuity of $F$ implies that for all $x\in[a,\infty]$ there exist a $b\in[a,\infty]$ such that $F(a,b)=x$, leading to $$F(a,x)=F(a,F(a,b))=F(F(a,a),b)=F(a,b)=x.$$ Thus $F(a,x)=\max\{a,x\}$ for all $x\in[0,\infty]$.
\end{proof}

\section{Weak pseudo-inverses of monotone functions}
This section first introduces a weak pseudo-inverse of a monotone function, and then discuss the properties of the weak pseudo-inverses. We also use the weak pseudo-inverse of a monotone function to construct a t-norm and a t-supconorm, respectively.
\begin{definition}\label{def3.10}
\emph{Let $a, b, m, n\in [-\infty, \infty]$ with $a<b, m<n$ and $t:[a,b]\rightarrow[m,n]$ be a monotone function. Then the function $t^{[-1]}:[m,n]\rightarrow[a,b]$ defined by
\begin{equation*}
t^{[-1]}(y)=\sup\{x\in [a,b]\mid (t(x)- y)(t(b)-t(a))\leq 0\}
\end{equation*}
is called a weak pseudo-inverse of the monotone function $t$.}
\end{definition}

 As an immediate consequence of Definition \ref{def3.10}, we get the following corollary.

\begin{corollary}\label{cor3.1}
 Let $a, b, m, n\in [-\infty, \infty]$ with $a<b, m<n$ and $t:[a,b]\rightarrow[m,n]$ be a monotone function
\renewcommand{\labelenumi}{(\roman{enumi})}
\begin{enumerate}
\item If $t$ is non-decreasing and non-constant, then for all $y\in [m,n]$ we obtain the simpler formula
 $$t^{[-1]}(y)=\sup\{x\in [a,b]\mid t(x)\leq y\}.$$
\item If $f$ is non-increasing and non-constant, then for all $y\in [m,n]$ we obtain the simpler formula
 $$t^{[-1]}(y)=\sup\{x\in [a,b]\mid t(x)\geq y\}.$$
\item If $t$ is a constant function, then for all $y\in [m,n]$ we have $t^{[-1]}(y)=b.$
\end{enumerate}
\end{corollary}

\begin{remark}\label{rem3.1}
\emph{Let $a, b, m, n\in [-\infty, \infty]$ with $a<b, m<n$ and $t:[a,b]\rightarrow[m,n]$ be a monotone function, and let $t^{[-1]}$ be its weak pseudo-inverse.
\renewcommand{\labelenumi}{(\roman{enumi})}
\begin{enumerate}
\item If $t$ is non-decreasing, then the function $t^{[-1]}$ is right continuous and non-decreasing, and for all $y\in[m,t(a))$ we get $t^{[-1]}(y)=a$, and for all $y\in(t(b),n]$ we have $t^{[-1]}(y)=b$.
\item If $t$ is non-increasing, then the function $t^{[-1]}$ is left continuous and non-increasing, and for all $y\in[m,t(b))$ we get $t^{[-1]}(y)=b$, and for all $y\in(t(a),n]$ we have $t^{[-1]}(y)=a$.
\end{enumerate}}
\end{remark}
\begin{example}\label{exp3.1}
\emph{Let the function $t:[0,1]\rightarrow [0,\infty]$ be defined by
\begin{equation*}
t(x)=\begin{cases}
x & \hbox{if }\ x\in[0,\frac{1}{4}],\\
\frac{1}{4} &  \hbox{if }\ x\in(\frac{1}{4},\frac{1}{2}],\\
x+ \frac{1}{2}&  \hbox{if }\ x\in(\frac{1}{2},\frac{3}{4}),\\
2 & \hbox{if }\ x\in[\frac{3}{4},\frac{7}{8}),\\
x+ \frac{5}{4}&  \hbox{if }\ x\in[\frac{7}{8},1].\\
\end{cases}
\end{equation*}
Then \begin{equation*}
t^{(-1)}(x)=\begin{cases}
x & \hbox{if }\ x\in[0,\frac{1}{4}],\\
\frac{1}{2} &  \hbox{if }\ x\in(\frac{1}{4},1],\\
x- \frac{1}{2}&  \hbox{if }\ x\in(1,\frac{5}{4}),\\
\frac{3}{4} & \hbox{if }\ x\in[\frac{5}{4},2],\\
\frac{7}{8} & \hbox{if }\ x\in(2,\frac{17}{8}),\\
x- \frac{5}{4}&  \hbox{if }\ x\in[\frac{17}{8},\frac{9}{4}),\\
1&  \hbox{if }\ x\in[\frac{9}{4},\infty].\\
\end{cases}
\end{equation*}
\begin{equation*}
t^{[-1]}(x)=\begin{cases}
x & \hbox{if }\ x\in[0,\frac{1}{4}),\\
\frac{1}{2} &  \hbox{if }\ x\in[\frac{1}{4},1],\\
x- \frac{1}{2}&  \hbox{if }\ x\in(1,\frac{5}{4}),\\
\frac{3}{4} & \hbox{if }\ x\in[\frac{5}{4},2),\\
\frac{7}{8} & \hbox{if }\ x\in[2,\frac{17}{8}),\\
x- \frac{5}{4}&  \hbox{if }\ x\in[\frac{17}{8},\frac{9}{4}),\\
1&  \hbox{if }\ x\in[\frac{9}{4},\infty].\\
\end{cases}
\end{equation*}
Obviously, $t^{(-1)}\leq t^{[-1]}$ and $t^*\leq t^{[-1]}$. If $x\in[0,\frac{1}{4})\cup (\frac{1}{2},\frac{3}{4})\cup[\frac{7}{8},1]$, then $t^{(-1)}(t(x))= t^{[-1]}(t(x))$. If $x\in[0,\frac{3}{4})\cup[\frac{7}{8},1]$, then $t^{[-1]}(t(x))=t^*(t(x))$.}
\end{example}
\begin{lemma}\label{lem3.1}
Let $t:[a,b]\rightarrow[m,n]$ be a non-decreasing (resp. non-increasing)
function and $t^{[-1]}$ be its weak pseudo-inverse.
Then for all $x\in[a,b]$,
$$t(t^{[-1]}(t(x)))\geq t(x) \mbox{ (resp. }t(t^{[-1]}(t(x)))\leq t(x)\mbox{)}.$$
\end{lemma}
\begin{proof}If $t$ is non-decreasing and $x\in[a,b]$ then
$$t^{[-1]}(t(x))=\sup\{y\in [a,b]\mid t(y)\leq t(x)\}\geq x,$$
thus,
$$t(t^{[-1]}(t(x)))=t(\sup\{y\in [a,b]\mid t(y)\leq t(x)\})\geq t(x).$$

The case that $t$ is non-increasing is completely analogous.
\end{proof}
\begin{lemma}\label{lem3.2}
Let $t:[a,b]\rightarrow[m,n]$ be a non-decreasing (resp. non-increasing)
function, let $t^{[-1]}$ be its weak pseudo-inverse and $x_{0}\in[a,b]$.
Then $t(t^{[-1]}(t(x_{0})))>t(x_{0})$ $(\mbox{resp. }t(t^{[-1]}(t(x_{0})))<t(x_{0}))$ if and only if there
exists a $\delta_{x_{0}}>0$ such that $t(x)=t(x_{0})$ for all
$x \in [x_{0},x_{0}+\delta_{x_{0}})$ and $t(x_{0})<t(x_{0} +\delta_{x_{0}})$ $(\mbox{resp. }t(x_{0})>t(x_{0} +\delta_{x_{0}}))$.
\end{lemma}
\begin{proof}
$(\Leftarrow)$. Let $t$ be a non-decreasing
function and $x_{0} \in[a,b)$. If there is a $\delta_{x_{0}}>0$ such
that $t(x)=t(x_{0})$ for all
$x \in [x_{0},x_{0}+\delta_{x_{0}})$ and $t(x_{0})<t(x_{0} +\delta_{x_{0}})$, then
$$t^{[-1]}(t(x_{0}))=\sup\{y\in [a,b]\mid t(y)\leq t(x_{0})\}= x_{0} +\delta_{x_{0}},$$
thus
$$t(t^{[-1]}(t(x_{0})))=t(\sup\{y\in [a,b]\mid t(y)\leq t(x_{0})\})=t(x_{0} +\delta_{x_{0}})>t(x_{0}).$$

$(\Rightarrow)$. Suppose that $t(t^{[-1]}(t(x_{0})))>t(x_{0})$.
Let $x_{0} \in[a,b)$ and $\alpha=t^{[-1]}(t(x_{0}))$.
Then $t(\alpha) > t(x_{0})$.
Using the monotonicity of $t$ we have $\alpha > x_{0}$.
Take $\delta_{x_{0}}=\alpha-x_{0}$. Then from the monotonicity of
$t$ and $\alpha=\sup\{y\in [a,b]\mid t(y)\leq t(x_{0})\}$, we immediately have
$t(x) = t(x_{0})$ for all
$x \in [x_{0},x_{0}+\delta_{x_{0}})$ and $t(x_{0})<t(x_{0} +\delta_{x_{0}})$.

The case that $t$ is non-increasing is completely analogous.
\end{proof}

 From Lemmas \ref{lem3.1} and \ref{lem3.2} we easily get the following theorem.
\begin{theorem}\label{th3.1}
Let $t:[a,b]\rightarrow[m,n]$ be a non-decreasing (resp. non-increasing)
function and $t^{[-1]}$ be its weak pseudo-inverse. Then the following are equivalent:
\renewcommand{\labelenumi}{(\roman{enumi})}
\begin{enumerate}
\item $\{x_{0}\in [a,b)\mid \mbox{ there is a } \delta_{x_{0}}>0 \mbox{ such that } t(x)=t(x_{0})
~\mbox{for all}~ x \in [x_{0},x_{0}+\delta_{x_{0}})~ \mbox{and}~ t(x_{0})<t(x_{0} +\delta_{x_{0}})\mbox{ (resp. }t(x_{0})>t(x_{0} +\delta_{x_{0}})\mbox{)}\}=\emptyset$;
\item $t(t^{[-1]}(t(x_{0})))=t(x_{0})$ for all $x_{0} \in [a,b]$.
\end{enumerate}
\end{theorem}

 From Definitions \ref{def2.3} and \ref{def3.10}, Lemma \ref{lem3.1} and
 Theorem \ref{th3.1}, we easily deduce the following properties of a weak pseudo-inverse of a monotone function.
\begin{proposition}\label{prop3.1}
Let $a, b, m, n\in [-\infty, \infty]$ with $a<b, m<n$ and $t:[a,b]\rightarrow[m,n]$ be a monotone function and $t^{[-1]}$ be its weak pseudo-inverse.
\renewcommand{\labelenumi}{(\roman{enumi})}
\begin{enumerate}
\item $t^{[-1]}$ coincides with $t^{(-1)}$ if and only if $t$ is strictly monotone. Moreover, $t^{[-1]}$ coincides with $t^{-1}$ if and only if $t$ is a bijection.
\item $t^{[-1]}$ is continuous  if and only if $t$ is strictly monotone on the set $t^{[-1]}([m,n))$.
\item $t^{[-1]}\circ t\geq id_{[a,b]}$.
\item If $t$ is either left continuous or strictly monotone then $t\circ t^{[-1]}\circ t=t$.
\item If $t$ is strictly monotone then so is $t^{[-1]}\mid_{Ran(t)}$. Further, we have
$$t\circ t^{[-1]}\mid_{Ran(t)}=id_{Ran(t)},~~ t^{[-1]}\circ t=id_{[a,b]}.$$
\item If $t$ is surjective then $t\circ t^{[-1]}=id_{[m,n]}$.
\item If both $\mu:[a,b]\rightarrow[a,b]$ and $\nu:[m,n]\rightarrow[m,n]$ are monotone bijections then
 $$(t\circ \mu)^{[-1]}=\mu^{-1}\circ t^{[-1]},~~(\nu\circ t)^{[-1]}= t^{[-1]}\circ \nu^{-1}.$$
\end{enumerate}
\end{proposition}

From Proposition \ref{prop3.1} (i) and Theorem \ref{th2.1}, we have the following corollary.
\begin{corollary}\label{cor3.2}
 Let $t :[0,1]\rightarrow [0,\infty]$ be a strictly decreasing function with $t(1) = 0$ such that
\begin{equation*}
t(x)+t(y)\in \emph{Ran}(t)\cup [t(0^{+}),\infty]
\end{equation*}
for all $(x,y)\in[0,1]^{2}$. Then the function $T: [0,1]^2\rightarrow [0,1]$ given by
\begin{equation*}
 T(x,y)= t^{[-1]}(t(x)+ t(y))
\end{equation*}
is a t-norm.
\end{corollary}

\begin{proposition}\label{prop3.2}
Let $t :[0,1]\rightarrow [0,\infty]$ be a left continuous non-decreasing function such that
\begin{equation}
\label{eq:(9)}
t(x)+t(y)\in \emph{Ran}(t)\cup [t(1^{-}),\infty]
\end{equation}
for all $(x,y)\in[0,1]^{2}$. Then the function $T: [0,1]^2\rightarrow [0,1]$ given by
\begin{equation*}
 T(x,y)= t^{[-1]}(t(x)+ t(y))
\end{equation*}
is a t-supconorm.
\end{proposition}
\begin{proof}Replacing $f^{(-1)}$ by $t^{[-1]}$, in completely analogous to the proof of Theorem 3.23 in \cite{EP2000} we can show the monotonicity, the associativity and the commutativity of $T$, respectively. On the other hand, $T(x,y)= t^{[-1]}(t(x)+ t(y))\geq t^{[-1]}(t(x))\geq x$ for all $x,y\in[0,1]$,
 analogously, $T(x,y)\geq y$. Thus $T(x,y)\geq \max\{x,y\}$. Therefore, by Definition \ref{de2.2} $T$ is a t-supconorm.
\end{proof}

Note that if $t :[0,1]\rightarrow [0,\infty]$ is a left continuous non-decreasing function but not strictly increasing and satisfies \eqref{eq:(9)} then one easily check that the function $T: [0,1]^2\rightarrow [0,1]$ given by $T(x,y)= t^{(-1)}(t(x)+ t(y))$ isn't a t-supconorm.

Generally, we can prove the following result through a analogous way to the proof of Proposition \ref{prop3.2}.
\begin{proposition}\label{prop3.3}
 Let $t :[0,1]\rightarrow [0,\infty]$ be a left continuous non-decreasing function and $F:[0,\infty]^2\rightarrow[0,\infty]$ be such that $([0,\infty],F,\leq)$
is a fully ordered Abel semigroup with $F(x,0) \geq x$ for all $x \in [0,\infty]$. If
\begin{equation}
F(t(x),t(y))\in \emph{Ran}(t)\cup [t(1^{-}),\infty]
\end{equation}
for all $(x,y)\in[0,1]^{2}$. Then the function $T: [0,1]^2\rightarrow [0,1]$ given by
\begin{equation}
 T(x,y)= t^{[-1]}F(t(x), t(y))
\end{equation}
is a t-supconorm.
\end{proposition}

In what follows, we consider what is a characterization of left continuous non-decreasing functions $t: [0,1]\rightarrow [0,\infty]$ such that the function $T: [0,1]^2\rightarrow [0,1]$ given by
\begin{equation}\label{eq:12}
T(x,y)= t^{[-1]}(F(t(x),t(y)))
\end{equation}
is associative, where $F:[0,\infty]^2\rightarrow[0,\infty]$ is an associative function and $t^{[-1]}:[0,\infty]\rightarrow[0,1]$ is the weak pseudo-inverse of $t$.

\section{The range of a left continuous non-decreasing function}

In this section we give a representation of the range of a left continuous non-decreasing function.

Let $t:[0,1]\rightarrow [0,\infty]$ be a function. We write $t(a^-)=\lim_{x \rightarrow a^{-}}t(x)$ for each $a\in (0, 1]$ and $t(a^+)=\lim_{x\rightarrow a^{+}}t(x)$ for each $a\in [0, 1)$. Define  $t(1^+)=\infty$ whenever $t$ is non-decreasing. Further, let
$$\mathcal{A}=\{M \mid \mbox{there is a left continuous non-decreasing function } t:[0,1]\rightarrow[0,\infty]\mbox{ such that }\mbox{Ran}(t)=M\}$$
and denoted by $A\setminus B=\{x\in A\mid x\notin B\}$ for two sets $A$ and $B$. Then in completely analogous to Lemma 3.1 of \cite{2025chen} we have the following lemma which presents the range of a left continuous non-decreasing function.

\begin{lemma}\label{lem4.1}
Let $t: [0,1]\rightarrow [0,\infty]$ be a left continuous non-decreasing function and $M\in \mathcal{A}$ with $M\neq[t(0), \infty]$. Then there exist a uniquely determined non-empty countable system $\mathcal{U}=\{[b_k, d_k] \subseteq [0,\infty]\mid k\in K\}$ of closed intervals of a positive length which satisfy that for all $[b_k, d_k],[b_l, d_l]\in \mathcal{U}$, $[b_k, d_k]\cap [b_l, d_l]=\emptyset$ or $[b_k, d_k]\cap [b_l, d_l]=\{d_k\}$ when $d_k\leq b_l$, and a uniquely determined non-empty countable set $\mathcal{V}=\{c_k\in [0,\infty]\mid k\in \overline{K}\}$ such that $[b_k, d_k] \cap \mathcal{V}=\{b_k\}$ or $[b_k, d_k] \cap \mathcal{V}=\{b_k,d_k\}$ for all $k\in K$ and
\begin{equation*}
M= \{c_k\in [0,\infty]\mid k\in \overline{K}\}\cup \left([t(0),\infty]\setminus \left(\bigcup_{k\in K}[b_k, d_k] \right)\right)
\end{equation*}
where $|K|\leq|\overline{K}|$.
\end{lemma}

\begin{definition}\label{def4.1}
\emph{Let $M\in \mathcal{A}$. A pair $(\mathcal{U},\mathcal{V})$ is said to be associated with $M\neq [t(0),\infty]$ if $\mathcal{U}=\{[b_k, d_k] \subseteq [0,\infty]\mid k\in K\}$ is a non-empty countable system of closed intervals of a positive length which satisfy that for all $[b_k, d_k],[b_l, d_l]\in \mathcal{U}$, $[b_k, d_k]\cap [b_l, d_l]=\emptyset$ or $[b_k, d_k]\cap [b_l, d_l]=\{d_k\}$ when $d_k\leq b_l$, and $\mathcal{V}=\{c_k\in [0,\infty]\mid k\in \overline{K}\}$ is a non-empty countable set such that $[b_k, d_k] \cap \mathcal{V}=\{b_k\}$ or $[b_k, d_k] \cap \mathcal{V}=\{b_k,d_k\}$ for all $k\in K$ and
\begin{equation*}
M= \{c_k\in [0,\infty]\mid k\in \overline{K}\}\cup \left([t(0),\infty]\setminus \left(\bigcup_{k\in K}[b_k, d_k] \right)\right).
\end{equation*}
A pair $(\mathcal{U},\mathcal{V})$ is said to be associated with $M=[t(0),\infty]$ if $\mathcal{U}=\{[\infty,\infty]\}$ and $\mathcal{V}=\{\infty\}$.}
\end{definition}

We briefly write $(\mathcal{U},\mathcal{V})=(\{[b_k, d_k]\mid k\in K\}, \{c_k \mid k\in \overline{K}\})$ instead of $(\mathcal{U},\mathcal{V})=(\{[b_k, d_k] \subseteq [0, \infty]\mid k\in K\}, \{c_k\in [0,\infty]\mid k\in \overline{K}\})$.

\begin{example}\label{exp4.1}\emph{
\renewcommand{\labelenumi}{(\roman{enumi})}
\begin{enumerate}
\item Let the function $t_1:[0,1]\rightarrow [0,\infty]$ be defined by\begin{equation*}
 t_1(x)=\begin{cases}
\frac{1}{2}x & \hbox{if }\ x\in[0,\frac{1}{2}],\\
\frac{1}{2} &  \hbox{if }\ x\in(\frac{1}{2},\frac{3}{4}],\\
x &  \hbox{if }\ x\in(\frac{3}{4},1].\\
\end{cases}
\end{equation*}
Then the pair $(\{[\frac{1}{4},\frac{1}{2}],[\frac{1}{2},\frac{3}{4}],[1,\infty]\},\{\frac{1}{4},\frac{1}{2},1\})$ is associated with $[0,\frac{1}{4}]\cup \{\frac{1}{2}\}\cup (\frac{3}{4},1]\in \mathcal{A}$.
\item Let the function $t_2:[0,1]\rightarrow [0,\infty]$ be defined by \begin{equation*}
 t_2(x)=\begin{cases}
1+x & \hbox{if }\ x\in[0,\frac{1}{5}],\\
\frac{6}{5} &  \hbox{if }\ x\in(\frac{1}{5},\frac{1}{4}],\\
\frac{3}{2} &  \hbox{if }\ x\in(\frac{1}{4},\frac{1}{2}],\\
2+x &  \hbox{if }\ x\in(\frac{1}{2},\frac{3}{4}],\\
\frac{1}{1-x} &  \hbox{if }\ x\in(\frac{3}{4},1),\\
\infty & \hbox{otherwise }.\\
\end{cases}
\end{equation*}
Then the pair $(\{[\frac{6}{5},\frac{3}{2}],[\frac{3}{2},\frac{5}{2}],[\frac{11}{4},4]\},\{\frac{6}{5},\frac{3}{2},\frac{11}{4}\})$ is associated with $[1,\frac{6}{5}]\cup \{\frac{3}{2}\}\cup [\frac{5}{2},\frac{11}{4}]\cup [4,\infty]\in \mathcal{A}$.
\end{enumerate}}
\end{example}

\section{Associativity of the function $T$ defined by Eq.(\ref{eq:12})}
This section shows necessary and sufficient conditions for the function $T$ defined by Eq.(\ref{eq:12}) being associative.

\subsection{An operation on $\mbox{Ran}(t)$ and its properties}

In this subsection we first define an operation $\otimes$ on $\mbox{Ran}(t)$ with $t$ a left continuous non-decreasing function, and then establish some necessary and sufficient conditions for the operation $\otimes$ being associative.

\begin{definition}\label{def5.2}
\emph{Let $M\in \mathcal{A}$. Define a function $G_{M}:[0,\infty]\rightarrow M $ by
\begin{equation*}
 G_{M}(x)=\min\{M\cap [\sup([0,x]\cap M), \inf([x,\infty]\cap M)]\}
\end{equation*}for all $\ x\in [0,\infty]$.}
\end{definition}

One can easily check the next proposition that describes the relationship between $M$ and $G_{M}$.

\begin{proposition}\label{prop5.1}
Let $M\in \mathcal{A}$ and $(\mathcal{U},\mathcal{V})=(\{[b_k, d_k]\mid k\in K\}, \{c_k \mid k\in \overline{K}\})$ be associated with $M$. Then for all $x\in[0,\infty]$ and $k\in K$,
\renewcommand{\labelenumi}{(\roman{enumi})}
\begin{enumerate}
\item If $x\in [0, t(0)]$ then $G_{M}(x)=t(0)$.
\item $G_{M}(x)=x$ if and only if $x\in M$.
\item If $x\notin M$ and $x> t(0)$ then $G_{M}(x)=b_k$ if and only if  $x\in [b_k, d_k]\setminus\{c_{k}\}$.
\item $G_{M}$ is a non-decreasing function.
\end{enumerate}
\end{proposition}
%\begin{proof}
% From Definition \ref{def5.2}, (i), (ii) and (iv) are immediate.
%
%(iii) Let $x\notin M$ and $x> t(0)$. Then there is a $k\in K$ such that $x\in [b_k, d_k]\setminus\{c_{k}\}$. Conversely, let $x> t(0)$ and $x\in [b_k, d_k]\setminus\{c_{k}\}$. Then $\sup([0,x]\cap M)=b_k$ and $\inf([x,1]\cap M)=d_k$. Thus from Definition \ref{def5.2}, $G_{M}(x)=b_{k}$.
%\end{proof}

\begin{example}\label{exp5.2}
\emph{In Example \ref{exp4.1},
\renewcommand{\labelenumi}{(\roman{enumi})}
\begin{enumerate}
\item \begin{equation*}
 G_{M}(x)=\begin{cases}
\frac{1}{4} & \hbox{if }\ x\in(\frac{1}{4},\frac{1}{2}),\\
\frac{1}{2} &  \hbox{if }\ x\in(\frac{1}{2},\frac{3}{4}],\\
x &  \mbox{otherwise}.\\
\end{cases}
\end{equation*}
\item \begin{equation*}
 G_{M}(x)=\begin{cases}
\frac{6}{5} & \hbox{if }\ x\in(\frac{6}{5},\frac{3}{2}),\\
\frac{3}{2} & \hbox{if }\ x\in(\frac{3}{2},\frac{5}{2}),\\
\frac{11}{4} &  \hbox{if }\ x\in(\frac{11}{4},4],\\
x &  \mbox{otherwise}.\\
\end{cases}
\end{equation*}
\end{enumerate}}
\end{example}

In what follows, we always suppose that $F:[0,\infty]^2\rightarrow [0,\infty]$ is an associative function. We need the following definition.

\begin{definition}\label{def5.3}
\emph{Let $M\in\mathcal{A}$ and $G_{M}$ be determined by $M$. Define an operation $\otimes:M^2\rightarrow M $ by
\begin{equation*}
x\otimes y=G_{M}(F(x,y)).
\end{equation*}}
\end{definition}
\begin{example}\label{exp5.3}
\emph{In Example \ref{exp5.2},
\renewcommand{\labelenumi}{(\roman{enumi})}
\begin{enumerate}
\item \begin{equation*}
 x\otimes y=\begin{cases}
\frac{1}{4} & \hbox{if }\ F(x,y)\in(\frac{1}{4},\frac{1}{2}),\\
\frac{1}{2} &  \hbox{if }\ F(x,y)\in(\frac{1}{2},\frac{3}{4}],\\
F(x,y) &  \hbox{otherwise}.\\
\end{cases}
\end{equation*}
\item \begin{equation*}
 x\otimes y=\begin{cases}
\frac{6}{5} & \hbox{if }\ F(x,y)\in(\frac{6}{5},\frac{3}{2}),\\
\frac{3}{2} &  \hbox{if }\ F(x,y)\in(\frac{3}{2},\frac{5}{2}),\\
\frac{11}{4} &  \hbox{if }\ F(x,y)\in(\frac{11}{4},4],\\
F(x,y) &  \hbox{otherwise}.\\
\end{cases}
\end{equation*}
\end{enumerate}}
\end{example}

\begin{proposition}\label{prop5.2}
Let $M\in \mathcal{A}$ and $(\mathcal{U},\mathcal{V})=(\{[b_k, d_k]\mid k\in K\}, \{c_k \mid k\in \overline{K}\})$ be associated with $M$. Then for all $x,y\in M$ and $k\in K$,
\renewcommand{\labelenumi}{(\roman{enumi})}
\begin{enumerate}
\item If $F(x,y)\in [0,t(0)]$,  then $x\otimes y=t(0)$.
\item $x\otimes y=F(x,y)$ if and only if $F(x,y)\in M$.
\item If $F(x,y)\notin M$ and $F(x,y)>t(0)$ then $x\otimes y=b_k$ if and only if  $F(x,y)\in [b_k, d_k]\setminus\{c_{k}\}$.
\item $\otimes$ is a non-decreasing function.
\end{enumerate}
\end{proposition}
\begin{proof}It is an immediate matter of Proposition \ref{prop5.1} and Definition \ref{def5.3}.
\end{proof}

\begin{proposition}\label{prop5.3}
Let $t:[0,1]\rightarrow [0,\infty] $ be a left continuous non-decreasing function with $\emph{Ran}(t)=M$ and $M\in \mathcal{A}$, $(\mathcal{U},\mathcal{V})=(\{[b_k, d_k]\mid k\in K\}, \{c_k \mid k\in \overline{K}\})$ be associated with $M$. Then $G_{M}(x)=t(t^{[-1]}(x))$ for all $x\in[0,\infty]$.
\end{proposition}
\begin{proof}If $x\in \mbox{Ran}(t)$ then, from Proposition \ref{prop3.1} (iv), we have $t(t^{[-1]}(x))=x$. Thus, in the case of $x\notin M$ and $x<t(0)$, by Remark \ref{rem3.1} we have $t^{[-1]}(x)=0$, hence $t(t^{[-1]}(x))=t(0)$; in the case of $x\notin M$ and $x>t(0)$, there is $k\in K$ such that $x\in [b_k, d_k]\setminus\{c_{k}\}$ with $b_k \in M$. Consequently,
\begin{eqnarray*}
t(t^{[-1]}(x))&=&t(\sup\{y\in [0,\infty]\mid t(y)\leq x\})\\
&=&\sup\{t(y)\in [0,\infty]\mid t(y)\leq x\}\\
&=&b_k.
\end{eqnarray*}
 Therefore, by Proposition \ref{prop5.1}, we get $G_{M}(x)=t(t^{[-1]}(x))$.
\end{proof}

Let $t:[0,1]\rightarrow [0,\infty]$ be a left continuous non-decreasing function with $\mbox{Ran}(t)=M$ and $(\mathcal{U},\mathcal{V})=(\{[b_k, d_k]\mid k\in K\}, \{c_k \mid k\in \overline{K}\})$ be associated with $M$.
Denote
$$\mathbb{H}=\{c\mid \mbox{ there are an } x_0\in[0,1]\mbox{ and }\varepsilon >0\mbox{ such that }t|_{[x_0,x_0+\varepsilon]}=c\},$$
$$\mathbb{G}=\{\sup \{x\in[0,1]\mid t(x)=y\}\mid y\in \mathbb{H}\},~\mathbb{W}=\{x\in [0,1]\mid t(x)\in M\setminus \mathbb{H} \},$$
$$\mathbb{D}=\mathbb{G}\cup \mathbb{W}.$$

In particular, $t^{[-1]}(x) \in \mathbb{D}$ for all $x\in [0,\infty]$, and we have the following definition.
\begin{definition}\label{def5.4}
\emph{Let $t:[0,1]\rightarrow [0,\infty]$ be a left continuous non-decreasing function.
Define a function $t^\star:\mathbb{D}\rightarrow [0,\infty]$ by
\begin{equation*}
 t^\star(x)=t(x)  \hbox{ for all }\ x\in \mathbb{D},
\end{equation*}
and a two-place function $F^{\star}:\mathbb{D}^2\rightarrow \mathbb{D}$ by
 \begin{equation*}
F^{\star}(x,y)=t^{[-1]}(F(t^\star(x),t^\star(y)))
\end{equation*}}
\end{definition}
for all $x,y\in \mathbb{D}$, respectively.

Then we immediately have the following remark.
\begin{remark}\label{rem5.1}
\emph{Let $t:[0,1]\rightarrow [0,\infty]$ be a left continuous non-decreasing function. Then
\renewcommand{\labelenumi}{(\roman{enumi})}
\begin{enumerate}
\item $t^\star$ is a strictly increasing function.
\item $t^{[-1]}( t^\star(x))=x$ for all $x\in \mathbb{D}$.
\item $t^\star ( t^{[-1]}(x))=x$ for all $x\in [0,\infty]$.
\end{enumerate}}
\end{remark}

In completely analogous to Proposition 4.3 of \cite{2025chen} we have the following lemma.
\begin{lemma}\label{lem5.4}
Let $t:[0,1]\rightarrow [0,\infty]$ be a left continuous non-decreasing function. Then
 \begin{equation*}
x\otimes y=t^\star(F^{\star}(t^{[-1]}(x),t^{[-1]}(y)))
\end{equation*}
for all $x,y\in M$ and
 \begin{equation*}
F^{\star}(x,y)=t^{[-1]}(t^\star(x)\otimes t^\star(y))
\end{equation*}
for all $x,y\in \mathbb{D}$.
\end{lemma}
%\begin{proof} From Definitions \ref{def5.3} and \ref{def5.4}, Remark \ref{rem5.1} and Proposition \ref{prop5.3}, for all $x,y\in M$ we get
%\begin{eqnarray*}
%x\otimes y&=&G_{M}(F(x,y))\\
%&=&t(t^{[-1]}(F(x,y)))\\
%&=&t^\star(t^{[-1]}(F(x,y)))\\
%&=&t^\star(t^{[-1]}(F(t^\star(t^{[-1]}(x)),t^\star(t^{[-1]}(y)))))\\
%&=&t^\star(F^{\star}(t^{[-1]}(x),t^{[-1]}(y))).
%\end{eqnarray*}
%Because of $x,y\in M$, there exist two elements $u,v\in \mathbb{D}$ such that $t^\star(u)=x$, $t^\star(v)=y$. Thus, from Remark \ref{rem5.1}, we have
%\begin{equation*}
%t^\star(u)\otimes t^\star(v)=x\otimes y=t^\star(F^{\star}(t^{[-1]}(x),t^{[-1]}(y)))=t^\star(F^{\star}(u,v)).
%\end{equation*}
%This follows that $t^{[-1]}(t^\star(u)\otimes t^\star(v))=t^{[-1]}(t^\star(F^{\star}(u,v))=F^{\star}(u,v)$.
%\end{proof}

Furthermore, by Lemma \ref{lem5.4} we have the following proposition.
\begin{proposition}\label{prop5.4}
Let $t:[0,1]\rightarrow [0,\infty]$ be a left continuous non-decreasing function. Then the following are equivalent:
\renewcommand{\labelenumi}{(\roman{enumi})}
\begin{enumerate}
\item $\otimes$ is associative.
\item $F^{\star}$ is associative.
\end{enumerate}
\end{proposition}
%\begin{proof}Let $\otimes$ be associative. Then, by Lemma \ref{lem5.4}, for all $x,y,z\in \mathbb{D}$ we have
%\begin{eqnarray*}
%F^{\star}(F^{\star}(x, y),z)&=&t^{[-1]}(t^{\star}(F^{\star}(x, y))\otimes t^{\star}(z))\\
%&=&t^{[-1]}(t^{\star}\circ t^{[-1]}(t^{\star}(x)\otimes t^{\star}(y))\otimes t^{\star}(z))\\
%&=&t^{[-1]}(t^{\star}(x)\otimes t^{\star}(y)\otimes t^{\star}(z))\\
%&=&t^{[-1]}(t^{\star}(x)\otimes t^{\star}\circ t^{[-1]}(t^{\star}(y)\otimes t^{\star}(z)))\\
%&=&t^{[-1]}(t^{\star}(x)\otimes t^{\star}(F^{\ast}(y, z)))\\
%&=&F^{\star}(x,F^{\star}(y,z)).
%\end{eqnarray*}
%Let $F^{\star}$ be associative. Then, by Lemma \ref{lem5.4}, for all $x,y,z\in M$ we have
%\begin{eqnarray*}
%(x\otimes y)\otimes z&=&t^\star(F^{\star}(t^{[-1]}(x\otimes y),t^{[-1]}(z)))\\
%&=&t^\star(F^{\star}(t^{[-1]}\circ t^\star(F^{\star}(t^{[-1]}(x),t^{[-1]}(y))),t^{[-1]}(z)))\\
%&=&t^\star(F^{\star}((F^{\star}(t^{[-1]}(x),t^{[-1]}(y))),t^{[-1]}(z)))\\
%&=&t^\star(F^{\star}(t^{[-1]}(x), F^{\star}(t^{[-1]}(y),t^{[-1]}(z)))\\
%&=&t^\star(F^{\star}(t^{[-1]}(x), t^{[-1]}\circ t^\star(F^{\star}(t^{[-1]}(y),t^{[-1]}(z))))\\
%&=&t^\star(F^{\star}(t^{[-1]}(x),t^{[-1]}(y\otimes z))\\
%&=&x\otimes (y\otimes z).
%\end{eqnarray*}
%\end{proof}

From Definition \ref{def5.4} we have the following lemma.
\begin{lemma}\label{lem5.2}
Let $t:[0,1]\rightarrow [0,\infty]$ be a left continuous non-decreasing function and $T:[0,1]^2\rightarrow [0,1]$ be the function defined by Eq.(\ref{eq:12}). Then, for each $x,y \in [0,1]$, there are two elements $m,n\in \mathbb{D}$ such that $t^\star(m)=t(x)$, $t^\star(n)=t(y)$ and $T(x,y)=F^{\star}(m,n)$. In particular, $T(x,y)=F^{\star}(x,y)$ for all $x,y\in \mathbb{D}$.
\end{lemma}
%\begin{proof}If $x,y\in \mathbb{D}$ then from Definition \ref{def5.4} we have $t^\star(x)=t(x)$ and $t^\star(y)=t(y)$. If $x\notin \mathbb{D}$, then $t(x)\in \mathbb{H}$. Let $m=\max \{s\in[0,1]\mid t(s)=t(x), t(x)\in \mathbb{H}\}$. Obviously, $m\in \mathbb{D}$ and $t^\star(m)=t(x)$. Analogously, if $y\notin \mathbb{D}$ then there is an $n\in \mathbb{D}$ such that $t^\star(n)=t(y)$.
%
%Therefore, from Definition \ref{def5.4} we have
%\begin{eqnarray*}
%T(x,y)&=&t^{[-1]}(F(t(x),t(y)))\\
%&=&t^{[-1]}(F(t^\star(m),t^\star(n)))\\
%&=&F^{\star}(m,n)
%\end{eqnarray*}for arbitrarily $x,y \in [0,1]$.
%\end{proof}

The next proposition describes the relation between $T$ and $F^{\star}$. One can easily prove it by Lemma \ref{lem5.2}.
\begin{proposition}\label{prop5.5}
Let $t:[0,1]\rightarrow [0,\infty]$ be a left continuous non-decreasing function and $T:[0,1]^2\rightarrow [0,1]$ be the function defined by Eq.(\ref{eq:12}).
Then the following are equivalent:
\renewcommand{\labelenumi}{(\roman{enumi})}
\begin{enumerate}
\item $T$ is associative.
\item $F^{\star}$ is associative.
\end{enumerate}
\end{proposition}
%\begin{proof}Suppose that $F^{\star}$ is associative. Let us prove that $T(T(x,y),z)=T(x,T(y,z))$ for all $x,y,z\in[0,1]$.
%Let $x,y,z\in [0,1]$. Then by Lemma \ref{lem5.2}, there exist $m,n,v\in \mathbb{D}$ such that $t^\star(m)=t(x)$, $t^\star(n)=t(y)$ and $t^\star(v)=t(z)$, respectively, and $T(x,y)=F^{\star}(m,n)$, $T(x,z)=F^{\star}(m,v)$ and $T(y,z)=F^{\star}(n,v)$, respectively. Therefore,
% \begin{eqnarray*}
%T(T(x,y),z)&=&T(F^{\star}(m,n),z)\\
%&=&F^{\star}(F^{\star}(m,n),v)\\
%&=&F^{\star}(m,F^{\star}(n,v))\\
%&=&F^{\star}(m,T(y,z))\\
%&=&T(x,T(y,z)).
%\end{eqnarray*}
%
%Conversely, if $T$ is associative, i,e., $T(T(x,y),z)=T(x,T(y,z)$ for all $x,y,z \in [0,1]$, then by Eq.(\ref{eq:12}), we have $$t^{[-1]}(F(t\circ t^{[-1]} (F(t(x),t(y))),t(z)))=t^{[-1]}(F(t(x),t\circ t^{[-1]} (F(t(x),t(y))))).$$
% So that $$t^{[-1]}(F(t^\star\circ t^{[-1]} (F(t^\star(x),t^\star(y))),t^*(z)))=t^{[-1]}(F(t^\star(x),t^\star\circ t^{[-1]} (F(t^\star(x),t^\star(y)))))$$
%for all $x,y,z\in \mathbb{D}$. Thus from Definition \ref{def5.4}, $F^{\star}(F^{\star}(x,y),z)=F^{\star}(x,F^{\star}(y,z))$, therefore, $F^{\star}$ is associative.
%\end{proof}

The following is an immediate consequence of Propositions \ref{prop5.4} and \ref{prop5.5}.
\begin{proposition}\label{prop5.6}
Let $t:[0,1]\rightarrow [0,\infty]$ be a left continuous non-decreasing function and $T:[0,1]^2\rightarrow [0,1]$ be the function defined by Eq.(\ref{eq:12}).
Then $T$ is associative if and only if $\otimes$ is associative.
\end{proposition}

\subsection{Associativity of the operation $\otimes$}
This subsection is devoted to exploring some necessary and sufficient conditions for the operation $\otimes$ being associative, which answer what properties of $M$ are equivalent to the associativity of $\otimes$.

Let $M\subseteq [0,\infty]$. Define $O(M) = \bigcup_{x,y\in M}(\min\{x,y\}, \max\{x,y\}]$ when $M\neq \emptyset$ (where $(x,x]=\emptyset$), and $O(M)=\emptyset$ when $M=\emptyset$.
Let $\emptyset\neq A,B\subseteq [0,\infty]$. Denote $F(A,B)=\{F(x,y)\mid x\in A, y\in B\}$ and $F(\emptyset,A)=\emptyset=F(A,\emptyset)$.

\begin{definition}\label{def6.1}
\emph{Let $M\in \mathcal{A}$ and $(\mathcal{U},\mathcal{V})=(\{[b_k, d_k]\mid k\in K\}, \{c_k \mid k\in \overline {K}\})$ be associated with $M$. For all $y\in M$ and $k,l\in K$, set $K^{\ast}=K\cup\{\tau\}$ where $\tau\notin K$.
Define
$$M_{k}^{y}=\{x\in M \mid F(x,y)\in [b_k, d_k]\setminus\{c_{k}\}\}, M_{y}^{k}=\{x\in M \mid F(y,x)\in [b_k, d_k]\setminus\{c_{k}\}\},$$
$$M_{\tau}^{y}=\{x\in M \mid F(x,y)< t(0)\}, M_{y}^{\tau}=\{x\in M \mid F(y,x)< t(0)\},$$
$$M^{y}=\{x\in M \mid F(x,y)\in M\}, M_{y}=\{x\in M \mid F(y,x)\in M\},$$
$$I_{k}^{y}=\{(x_1,x_2)\in M_{k}^{y}\times M_{y} \mid F(F(x_1,y),x_2)\neq F(c_k,x_2), ~ M_{k}^{y}\neq\emptyset,  M_{y}\neq\emptyset\},$$ 
			$$I_{\tau}^{y}=\{(x_1,x_2)\in M_{\tau}^{y}\times M_{y} \mid F(F(x_1,y),x_2)\neq F(t(0),x_2), ~ M_{\tau}^{y}\neq\emptyset,  M_{y}\neq\emptyset\},$$ 
		$$I_{y}^{k}=\{(x_1,x_2)\in M^{y} \times  M_{y}^{k}\mid F(x_1,F(y,x_2))\neq F(x_1, c_k), ~ M_{y}^{k}\neq\emptyset,  M^{y}\neq\emptyset\},$$ 
				$$I_{y}^{\tau}=\{(x_1,x_2)\in M^{y} \times  M_{y}^{\tau}\mid F(x_1,F(y,x_2))\neq F(x_1, t(0)), ~ M_{y}^{\tau}\neq\emptyset,  M^{y}\neq\emptyset\},$$ 
		$$I_{k,l}^{y}=\{(x_1,x_2)\in M_{k}^{y}\times M_{l}^{y} \mid F(c_k,x_2)\neq F(x_1,c_l), ~ M_{k}^{y}\neq\emptyset,  M_{l}^{y}\neq\emptyset\},$$
				$$I_{k,\tau}^{y}=\{(x_1,x_2)\in M_{k}^{y}\times M_{\tau}^{y} \mid F(c_k,x_2)\neq F(x_1,t(0)), ~ M_{k}^{y}\neq\emptyset,  M_{\tau}^{y}\neq\emptyset\},$$
						$$I_{\tau,l}^{y}=\{(x_1,x_2)\in M_{\tau}^{y}\times M_{l}^{y} \mid F(t(0),x_2)\neq F(x_1,c_l), ~ M_{\tau}^{y}\neq\emptyset,  M_{l}^{y}\neq\emptyset\},$$
				$$I_{\tau,\tau}^{y}=\{(x_1,x_2)\in M_{\tau}^{y}\times M_{\tau}^{y} \mid F(t(0),x_2)\neq F(x_1,t(0)), ~ M_{\tau}^{y}\neq\emptyset,  M_{\tau}^{y}\neq\emptyset\},$$
  $$H_{k}^{y}=O(\{b_k\}\cup F(M_{k}^{y},y)), H_{\tau}^{y}=O(\{t(0)\}\cup F(M_{\tau}^{y},y)),$$
$$H_{y}^{k}=O(\{b_k\}\cup F(y,M_{y}^{k})), H_{y}^{\tau}=O(\{t(0)\}\cup F(y,M_{y}^{\tau})),$$
\begin{equation*}
J_{k,l}^{y}=\begin{cases}
O(F(M_{k}^{y},b_l)\cup F(b_k,M_{y}^{l})), & \mbox{if } M_{k}^{y}\neq\emptyset,  M_{y}^{l}\neq\emptyset,\\
\emptyset, & \hbox{otherwise.}
\end{cases}
\end{equation*}
\begin{equation*}
J_{\tau,l}^{y}=\begin{cases}
O(F(M_{\tau}^{y},b_l)\cup F(t(0),M_{y}^{l})), & \mbox{if } M_{\tau}^{y}\neq\emptyset,  M_{y}^{l}\neq\emptyset,\\
\emptyset, & \hbox{otherwise.}
\end{cases}
\end{equation*}
\begin{equation*}
J_{k,\tau}^{y}=\begin{cases}
O(F(M_{k}^{y},t(0))\cup F(b_k,M_{y}^{\tau})), & \mbox{if } M_{k}^{y}\neq\emptyset,  M_{y}^{\tau}\neq\emptyset,\\
\emptyset, & \hbox{otherwise.}
\end{cases}
\end{equation*}
\begin{equation*}
J_{\tau,\tau}^{y}=\begin{cases}
O(F(M_{\tau}^{y},t(0))\cup F(t(0),M_{y}^{\tau})), & \mbox{if } M_{\tau}^{y}\neq\emptyset,  M_{y}^{\tau}\neq\emptyset,\\
\emptyset, & \hbox{otherwise.}
\end{cases}
\end{equation*}
Put
$\mathfrak{T_{1}}(M)=\bigcup_{y\in M}\bigcup_{k\in K^{\ast}}F(H_{k}^{y},M_{y}), \mathfrak{T_{2}}(M)=\bigcup_{y\in M}\bigcup_{k\in K^{\ast}}F(M^{y},H_{y}^{k}),  \mathfrak{T_{3}}(M)=\bigcup_{y\in M}\bigcup_{k,l\in K^{\ast}}J_{k,l}^{y}$
and $\mathfrak{T}(M)=\mathfrak{T_{1}}(M)\cup\mathfrak{T_{2}}(M)\cup\mathfrak{T_{3}}(M).$}
\end{definition}

In the rest of this section, we further suppose that $F:[0,\infty]^2\rightarrow [0,\infty]$ is a monotone and associative function with neutral element in $[0,\infty]$. Then similar to Lemmas 5.3 and 5.4 in \cite{YM2024}, respectively, by Proposition \ref{prop5.1} we have the following two lemmas.
\begin{lemma}\label{lem6.1}
Let $M\in \mathcal{A}$ and $(\mathcal{U},\mathcal{V})=(\{[b_k, d_k]\mid k\in K\}, \{c_k \mid k\in \overline{K}\})$ be associated with $M$. Then for any $x,y\in[0,\infty]$, $G_M(x)= G_M(y)$ if and only if $(\min\{x,y\},\max\{x,y\}]\cap (M\setminus t(0))=\emptyset$.
\end{lemma}
%\begin{proof}Let $G_M(x)= G_M(y)$. If $(\min\{x,y\},\max\{x,y\}]\cap (M\setminus t(0))\neq\emptyset$, then $x\neq y$, without loss of generality, say $x<y$, thus there is an $a \in(x,y]\cap (M\setminus t(0))$.
%By Proposition \ref{prop5.1}, $G_M(a)=a>t(0)$. If $x\in M$, then by Proposition \ref{prop5.1} $G_M(x)=x$, thus $G_M(x)< G_M(a)$. If $x\notin M$, then $x\in [0,t(0))$ or there is a $k\in K$ such that $x\in [b_k,d_k]\setminus\{c_{k}\}$. In the case $x\in [0,t(0))$, by Proposition \ref{prop5.1} we have $G_M(x)=t(0)< G_M(a)$. In the case $x\in [b_k,d_k]\setminus\{c_{k}\}$, by Proposition \ref{prop5.1} we have $G_M(x)=b_k$, thus $G_M(x)<x<a=G_M(a)$. In summary, we always have $G_M(x)< G_M(a)\leq G_M(y)$,
% a contradiction.
%
%Conversely, if $(\min\{x,y\},\max\{x,y\}]\cap (M\setminus t(0))=\emptyset$, without loss of generality, say $x< y$, then $(x,y]\subseteq (0,t(0)]$ or there is a $k\in K$ such that $(x,y]\subseteq [b_k,d_k]\setminus\{c_{k}\}$. If $(x,y]\subseteq (0,t(0)]$ then, from Proposition \ref{prop5.1},$G_M(x)=0=G_M(y)$. If $(x,y]\subseteq [b_k,d_k]\setminus\{c_{k}\}$ then, from Proposition \ref{prop5.1}, $G_M(x)=b_k=G_M(y)$.
%\end{proof}

\begin{lemma}\label{lem6.2}
Let $M\in \mathcal{A}$ and $(\mathcal{U},\mathcal{V})=(\{[b_k, d_k]\mid k\in K\}, \{c_k \mid k\in \overline{K}\})$ be associated with $M$. Let $M_1, M_2\subseteq [0,\infty]$ be two non-empty sets and $c\in[0,\infty]$. If there exist $u\in M_1$ and $v\in M_2$ such that $F(u,c)\neq F(v,c)$ and $F(c,u)\neq F(c,v)$,
Then\\
(1) $F(O(M_1\cup M_2),c)\cap (M\setminus \{t(0)\})\neq\emptyset$ if and only if there exist $x\in M_1$ and $y\in M_2$ such that $$(\min\{F(x,c),F(y,c)\},\max\{F(x,c),F(y,c)\}]\cap (M\setminus \{t(0)\}) \neq\emptyset.$$
(2) $F(c,O(M_1\cup M_2))\cap (M\setminus \{t(0)\})\neq\emptyset$ if and only if there exist $x\in M_1$ and $y\in M_2$ such that $$(\min\{F(c,x),F(c,y)\},\max\{F(c,x),F(c,y)\}]\cap (M\setminus \{t(0)\})\neq\emptyset.$$
\end{lemma}
%\begin{proof} In a completely analogous to the proof of Lemma 5.3 in \cite{YM2024}.
 %(1) We only provide the proof of statement (1) when $F$ is non-decreasing, the statement being analogous when $F$ is non-increasing.
%Suppose that $F(O(M_1\cup M_2),c)\cap (M\setminus \{t(0)\})\neq\emptyset$. Then there exist two elements $s,t\in M_1\cup M_2$ with $s<t$ such that $F((s,t],c)\cap (M\setminus \{t(0)\})\neq\emptyset$.
%
%The following proof is split into three cases.
%\renewcommand{\labelenumi}{(\roman{enumi})}
%\begin{enumerate}
%\item The case when one of $\{s,t\}$ is contained in $M_1$ and the other is contained in $M_2$ is obvious.
%
%\item The case that $s,t\in M_1$ with $s<t$. Choose $w\in M_2$. There are three subcases.
%\renewcommand{\labelenumi}{(\roman{enumi})}
%\begin{enumerate}
%\item If $w\leq s$ then put $x=t$ and $y =w$.
%\item If $w\geq t$ then put $x=s$ and $y =w$.
%\item If $s<w<t$, then $F((s,w],c)\cap (M\setminus \{t(0)\})\neq\emptyset$ or $F((w,t),c)\cap (M\setminus \{t(0)\})\neq\emptyset$. If $F((s,w],c)\cap (M\setminus \{t(0)\})\neq\emptyset$ then put $x=s$ and $y =w$. If $F([w,t),c)\cap (M\setminus \{t(0)\})\neq\emptyset$ then put $x=t$ and $y =w$.
%\end{enumerate}
%\item The case that $s,t\in M_2$ with $s<t$ is completely analogous to (ii).
%\end{enumerate}
%The converse implication is obvious.
%
%(2) The proof is analogy to the proof of statement (1).
%\end{proof}

Let $M\in \mathcal{A}$, $(\mathcal{U},\mathcal{V})=(\{[b_k, d_k]\mid k\in K\}, \{c_k \mid k\in \overline{K}\})$ be associated with $M$. For all  $k,l\in K^*, y\in M$, write

$(C_1)$ either $I_{k}^{y} = \emptyset$ or $F(H_{k}^{y}, M_{y}) \cap (M\setminus \{t(0)\}) = \emptyset$;

$(C_2)$ either $I_{y}^{k} = \emptyset$ or $F(M^{y},H_{y}^{k}) \cap (M\setminus \{t(0)\}) = \emptyset$;

$(C_3)$ either $I_{k,l}^{y} = \emptyset$ or $J_{k,l}^{y} \cap (M\setminus \{t(0)\}) = \emptyset$.

Conditions $(C_1)$, $(C_2)$ and $(C_3)$ are called an $F$-condition of $M$.

The following proposition characterizes what properties of $M$ are equivalent to the associativity of $\otimes$.
\begin{proposition}\label{prop5.7}
Let $M\in \mathcal{A}$ and $(\mathcal{U},\mathcal{V})=(\{[b_k, d_k]\mid k\in K\}, \{c_k \mid k\in \overline{K}\})$ be associated with $M$. Then  the operation $\otimes$ on $M$ is associative if and only if the $F$-condition of $M$ holds.
\end{proposition}
\begin{proof} Let $(\mathcal{U},\mathcal{V})=(\{[b_k, d_k]\mid k\in K\}, \{c_k \mid k\in \overline{K}\})$ be associated with $M\in \mathcal{A}$. We prove that the operation $\otimes$ on $M$ is not associative if and only if the $F$-condition of $M$ does not hold.

 Suppose that the operation $\otimes$ is not associative, i.e., there exist three elements $x,y,z\in M$ such that $(x\otimes y)\otimes z\neq x\otimes(y\otimes z)$. Then we claim that $F(x,y)\notin M$ or $F(y,z)\notin M$. Otherwise, from Proposition \ref{prop5.2}, $F(x,y)\in M$ and $F(y,z)\in M$ would imply $(x\otimes y)\otimes z=G_M(F(F(x,y),z))= G_M(F(x,F(y,z)))= x\otimes(y\otimes z)$, a contradiction. We consider three cases as below.

(i) Let $F(x,y)\notin M$ and $F(y,z)\in M$. Then $y\otimes z=F(y,z)$ and either $F(x,y)\in[0,t(0))$ or there exists a $k\in K$ such that $F(x,y)\in [b_k, d_k]\setminus\{c_{k}\}$. If $F(x,y)\in[0,t(0))$ then, by Proposition \ref{prop5.2}, $x\otimes y=t(0)$. It follows from Definition \ref{def5.3} that $G_M(F(t(0),z))=(x\otimes y)\otimes z\neq x\otimes(y\otimes z)=G_M(F(x,F(y,z)))$. On the other hand, by the associativity of $F$, we have $G_M(F(x,F(y,z)))=G_M(F(F(x,y),z))$. Thus $G_M(F(t(0),z))\neq G_M(F(F(x,y),z))$. So, $I_{\tau}^{y}\neq \emptyset$. Therefore, by Lemma \ref{lem6.1}, $$(\min\{F(t(0),z),F(F(x,y),z)\},\max\{F(t(0),z),F(F(x,y),z)\}]\cap (M\setminus \{t(0)\})\neq\emptyset .$$ Obviously, $H_{\tau}^{y}=O(\{t(0)\}\cup F(M_{k}^{y},y))$, $z\in M_{y}$ and $x\in M_{k}^{y}$, so that $$(\min\{F(t(0),z),F(F(x,y),z)\},\max\{F(t(0),z),F(F(x,y),z)\}]\subseteq F(H_{\tau}^{y},M_{y}),$$ which implies $F(H_{\tau}^{y},M_{y})\cap (M\setminus \{t(0)\}) \neq\emptyset$.

 If $F(x,y)\in [b_k, d_k]\setminus\{c_{k}\}$, then $x\otimes y=b_k$. It follows from Definition \ref{def5.2} that $G_M(F(b_k,z))=(x\otimes y)\otimes z\neq x\otimes(y\otimes z)=G_M(F(x,F(y,z)))$. On the other hand, by the associativity of $F$, we have $G_M(F(x,F(y,z)))=G_M(F(F(x,y),z))$. Thus $G_M(F(b_k,z))\neq G_M(F(F(x,y),z))$. So, $I_{k}^{y}\neq \emptyset$. Therefore, by Lemma \ref{lem6.1}, $$(\min\{F(b_k,z),F(F(x,y),z)\},\max\{F(b_k,z),F(F(x,y),z)\}]\cap (M\setminus \{t(0)\})\neq\emptyset. $$ Obviously, $H_{k}^{y}=O(\{b_k\}\cup F(M_{k}^{y},y))$, $z\in M_{y}$, and $x\in M_{k}^{y}$. So that $$(\min\{F(b_k,z),F(F(x,y),z)\},\max\{F(b_k,z),F(F(x,y),z)\}]\subseteq F(H_{k}^{y},M_{y}),$$ which implies $F(H_{k}^{y},M_{y})\cap (M\setminus \{t(0)\}) \neq\emptyset$.

(ii) Let $F(x,y)\in M$ and $F(y,z)\notin M$. In completely analogous to (i), $F(M^{y},H_{y}^{k})\cap (M\setminus \{t(0)\}) \neq\emptyset$ where $k\in K^{\ast}$.

(iii) Let $F(x,y)\notin M$ and $F(y,z)\notin M$. Then from $F(x,y)\notin M$, we have $F(x,y)\in[0,t(0))$ or there exists a $k\in K$ such that $F(x,y)\in [b_k, d_k]\setminus\{c_{k}\}$. In the case $F(x,y)\in[0,t(0))$, from Proposition \ref{prop5.2} we have $x\otimes y=t(0)$. In the case $F(x,y)\in [b_k, d_k]\setminus\{c_{k}\}$, we have $x\otimes y=b_k$. Hence $(x\otimes y)\otimes z=G_M(F(t(0),z))$ or $(x\otimes y)\otimes z=G_M(F(b_k,z))$. From $F(y,z)\notin M$, we have $F(y,z)\in[0,t(0))$ or there exists an $l\in K$ such that $F(y,z)\in [b_l,d_l]\setminus\{c_{l}\}$. If $F(y,z)\in[0,t(0))$ then, from Proposition \ref{prop5.2}, $y\otimes z=t(0)$. If $F(y,z)\in [b_l,d_l]\setminus\{c_{l}\}$ then $y\otimes z=b_l$. Thus $x\otimes(y\otimes z)=G_M(F(x,t(0)))$ or $x\otimes(y\otimes z)=G_M(F(x,b_l))$.

Since $(x\otimes y)\otimes z\neq x\otimes(y\otimes z)$, obviously, $I_{k,l}^{y} \neq \emptyset$ for all $k,l \in K^*$. By Lemma \ref{lem6.1} there are four cases as follows.

Case (1). $(\min\{F(t(0),z),F(x,t(0))\}, \max\{F(t(0),z),F(x,t(0))\}]\cap (M\setminus \{t(0)\})\neq \emptyset$. Obviously, $x\in M_{k}^{y}$ and $F(x,t(0))\in F(M_{k}^{y},t(0))$. Similarly, $z\in M_{y}^{l}$ and $F(t(0),z)\in F(t(0),M_{y}^{l})$. Therefore, $$(\min\{F(t(0),z),F(x,t(0))\},\max\{F(t(0),z),F(x,t(0))\}]\subseteq J_{\tau,\tau}^{y}.$$ This follows $J_{\tau,\tau}^{y}\cap (M\setminus \{t(0)\})\neq \emptyset$.

Case (2). $(\min\{F(t(0),z),F(x,b_l)\}, \max\{F(t(0),z),F(x,b_l)\}]\cap (M\setminus \{t(0)\})\neq \emptyset$. Obviously, $x\in M_{k}^{y}$ and $F(x,b_l)\in F(M_{k}^{y},b_l)$. Similarly, $z\in M_{y}^{l}$ and $F(t(0),z)\in F(t(0),M_{y}^{l})$. Therefore, $$(\min\{F(t(0),z),F(x,b_l)\},\max\{F(t(0),z),F(x,b_l)\}]\subseteq J_{\tau,l}^{y}.$$ This follows $J_{\tau,l}^{y}\cap (M\setminus \{t(0)\})\neq \emptyset$.

Case (3). $(\min\{F(b_k,z),F(x,t(0))\}, \max\{F(b_k,z),F(x,t(0))\}]\cap (M\setminus \{t(0)\})\neq \emptyset$. Obviously, $x\in M_{k}^{y}$ and $F(x,t(0))\in F(M_{k}^{y},t(0))$. Similarly, $z\in M_{y}^{l}$ and $F(b_k,z)\in F(b_k,M_{y}^{l})$. Therefore, $$(\min\{F(b_k,z),F(x,t(0))\},\max\{F(b_k,z),F(x,t(0))\}]\subseteq J_{k,\tau}^{y}.$$ This follows $J_{k,\tau}^{y}\cap (M\setminus \{t(0)\})\neq \emptyset$.

Case (4). $(\min\{F(b_k,z),F(x,b_l)\}, \max\{F(b_k,z),F(x,b_l)\}]\cap (M\setminus \{t(0)\})\neq \emptyset$. Obviously, $x\in M_{k}^{y}$ and $F(x,b_l)\in F(M_{k}^{y},b_l)$. Similarly, $z\in M_{y}^{l}$ and $F(b_k,z)\in F(b_k,M_{y}^{l})$. Therefore, $$(\min\{F(b_k,z),F(x,b_l)\},\max\{F(b_k,z),F(x,b_l)\}]\subseteq J_{k,l}^{y}.$$ This follows $J_{k,l}^{y}\cap (M\setminus \{t(0)\})\neq \emptyset$.

Finally, (i), (ii) and (iii) mean that the $F$-condition of $M$ does not hold.

Conversely, suppose the $F$-condition of $M$ does not hold. Then there exist a $y\in M$ and two elements $k,l\in K^{\ast}$ such that $I_{k}^{y}\neq \emptyset$ and $F(H_{k}^{y},M_{y})\cap (M\setminus \{t(0)\}) \neq\emptyset$, or $I_{y}^{k}\neq \emptyset$ and $F(M^{y},H_{y}^{k})\cap (M\setminus \{t(0)\}) \neq\emptyset$, or $I_{k,l}^{y}\neq \emptyset$ and  $J_{k,l}^{y}\cap (M\setminus \{t(0)\})\neq \emptyset$. We distinguish three cases as follows.

(i) Let $I_{k}^{y}\neq \emptyset$ and $F(H_{k}^{y},M_{y})\cap (M\setminus \{t(0)\}) \neq\emptyset$. If $k\in K^{\ast}$, then there exists a $z\in M_{y}$ such that $F(H_{k}^{y},z)\cap (M\setminus \{t(0)\}) \neq\emptyset$. Thus by the definition of $H_{k}^{y}$, $F(M_{k}^{y},y)\neq\emptyset$. Because of $I_{k}^{y}\neq \emptyset$, applying Lemma \ref{lem6.2}, there exist $u\in\{t(0),b_k\}$ and $v\in F(M_{k}^{y},y)$ such that
$$(\min\{F(u,z),F(v,z)\},\max\{F(u,z),F(v,z)\}]\cap (M\setminus \{t(0)\})\neq\emptyset.$$
Because of $v\in F(M_{k}^{y},y)$, there exists an $x\in M_{k}^{y}$ such that $F(x,y)=v$. Therefore, there exist two elements $u\in\{t(0),b_k\}$ and $x\in M_{k}^{y}$ such that
$$(\min\{F(u,z),F(F(x,y),z)\},\max\{F(u,z),F(F(x,y),z)\}]\cap(M\setminus \{t(0)\})\neq\emptyset.$$
Consequently, from Lemma \ref{lem6.1} we have $G_M(F(u,z))\neq G_M(F(F(x,y),z))$.
On the other hand, from $z\in M_{y}$, we have $F(y,z)\in M$. This follows $y\otimes z=F(y,z)$. From $x\in M_{k}^{y}$ we have $F(x,y)\in [0,t(0))$ or $F(x,y)\in [b_k,d_k]\setminus\{c_{k}\}$.
 If $F(x,y)\in [0,t(0))$ then $x\otimes y=t(0)$. Therefore, $(x\otimes y)\otimes z=G_M(F(t(0),z))\neq G_M(F(F(x,y),z))=G_M(F(x,F(y,z)))= x\otimes(y\otimes z)$.
If $F(x,y)\in [b_k,d_k]\setminus\{c_{k}\}$ then $x\otimes y=b_k$. Therefore, $(x\otimes y)\otimes z=G_M(F(b_k,z))\neq G_M(F(F(x,y),z))=G_M(F(x,F(y,z)))= x\otimes(y\otimes z)$.

(ii) Let $I_{y}^{k}\neq \emptyset$ and $F(M^{y},H_{y}^{k})\cap (M\setminus \{t(0)\}) \neq\emptyset$. Then in complete analogy to (i), $(x\otimes y)\otimes z\neq x\otimes(y\otimes z)$.

(iii) Let $I_{k,l}^{y}\neq \emptyset$ and  $J_{k,l}^{y}\cap (M\setminus \{t(0)\})\neq \emptyset$. Then $J_{k,l}^{y}\neq \emptyset$. Thus by the definition of $J_{k,l}^{y}$, we have $F(O(F(M_{k}^{y},a)\cup F(b,M_{y}^{l})),e)\cap (M\setminus \{t(0)\})\neq \emptyset$ where $a\in \{t(0),b_l\}$, $b\in \{t(0),b_k\}$ and $e$ is a neutral element of $F$, which means $F(M_{k}^{y},a)\neq\emptyset$ and $F(b,M_{y}^{l})\neq\emptyset$. Because of $I_{k,l}^{y}\neq \emptyset$, applying Lemma \ref{lem6.2}, there exist two elements $u\in F(b,M_{y}^{l})$ and $v\in F(M_{k}^{y},a)$ such that $$(\min\{F(u,e),F(v,e)\},\max\{F(u,e),F(v,e)\}]\cap (M\setminus \{t(0)\})\neq\emptyset.$$
Because $u\in F(b,M_{y}^{l})$ and $v\in F(M_{k}^{y},a)$, there exist a $z\in M_{y}^{l}$ and an $x\in M_{k}^{y}$ such that $u=F(b,z)$, $v=F(x,a)$. Therefore, there exist an $x\in M_{k}^{y}$ and a $z\in M_{y}^{l}$ such that
$$(\min\{F(b,z),F(x,a)\},\max\{F(b,z),F(x,a)\}]\cap (M\setminus \{t(0)\})\neq\emptyset$$
since $e$ is a neutral element of $F$. Further, by Lemma \ref{lem6.1} we have $G_M(F(b,z))\neq G_M(F(x,a))$.

On the other hand, from $x\in M_{k}^{y}$ we have  $F(x,y)\in [0,t(0))$ or $F(x,y)\in [b_k,d_k]\setminus\{c_{k}\}$. Thus $x\otimes y=t(0)$ or $x\otimes y=b_k$.
From $z\in M_{y}^{l}$ we have $F(y,z)\in[0,t(0))$ or $F(y,z)\in [b_l,d_l]\setminus\{c_{l}\}$. Thus $y\otimes z=t(0)$ or $y\otimes z=b_l$.

Therefore, $(x\otimes y)\otimes z=G_M(F(b,z))\neq G_M(F(x,a))= x\otimes(y\otimes z)$.
\end{proof}

Further, from Propositions \ref{prop5.6} and \ref{prop5.7} we have the following theorem.

\begin{theorem}\label{th6.2}
Let $t:[0,1]\rightarrow [0,\infty]$ be a left continuous non-decreasing function and $T:[0,1]^2\rightarrow [0,1]$ be the function defined by Eq.(\ref{eq:12}). Then the function $T$ is associative if and only if the $F$-condition of $M$ holds.
\end{theorem}

In particular, if $F$ is strictly monotone, then for all  $k,l\in K^*, y\in M$, $I_{k}^{y} =I_{y}^{k}=I_{k,l}^{y}=\emptyset$. Hence, we have the following corollary.
\begin{corollary}\label{co02.3.001}
Let $t:[0,1]\rightarrow [0,\infty]$ be a left continuous non-decreasing function with $\emph{Ran}(t)=M$ and $T:[0,1]^2\rightarrow [0,1]$  be a function defined by Eq.(\ref{eq:12}). If $F$ is strictly monotone, then the function $T$ is associative if and only if $\mathfrak{T}(M)\cap (M\setminus t(0))=\emptyset$.
\end{corollary}

From Definition \ref{def6.1}, we have $0\notin \mathfrak{T}(M)$. Thus if $t(0)=0$, then $\mathfrak{T}(M)\cap (M\setminus \{t(0)\})=\emptyset$ if and only if $\mathfrak{T}(M)\cap M =\emptyset$. Therefore, we have the following corollary.
\begin{corollary}\label{coro6.4}
Let $t:[0,1]\rightarrow [0,\infty]$ be a left continuous non-decreasing function with $t(0)=0$ and $\emph{Ran}(t)=M$. If $F$ is strictly monotone, then the function $T$ given by Eq.(\ref{eq:12}) is associative if and only if $\mathfrak{T}(M)\cap M=\emptyset$.
\end{corollary}

From Proposition \ref{prop3.3}, we have the following corollary.
\begin{corollary}\label{coro6.5}
Let $t:[0,1]\rightarrow [0,\infty]$ be a left continuous non-decreasing function and $T:[0,1]^2\rightarrow [0,1]$ be the function defined by Eq.(\ref{eq:12}). If $F:[0,\infty]^2\rightarrow[0,\infty]$ is a function such that $([0,\infty],F,\leq)$ is a fully ordered Abel semigroup with neutral element in $[0,\infty]$ and $F(x,0) \geq x$ for all $x \in [0,\infty]$, then $T$ is a t-supconorm if and only if the $F$-condition of $M$ holds.
\end{corollary}

%In particular, we have the following remark.
%\begin{remark}\label{rem3.3} \emph{Let $t:[0,1]\rightarrow [0,\infty]$ be a left continuous non-decreasing function and $T:[0,1]^2\rightarrow [0,1]$ be a function defined by Eq.(\ref{eq:12}).}
%\begin{enumerate}\renewcommand{\labelenumi}{(\roman{enumi})}
%\item \emph{ If $0$ is a neutral element of $F$, then $T(x,0)=t^{[-1]}(F(t(x),t(0)))\geq t^{[-1]}(t(x))\geq x$ for all $x\in [0,1]$. So, $0$ isn't necessary a neutral element of $T$. Therefore, if $\mathfrak{T}(M)\cap (M\setminus \{t(0)\})=\emptyset$ then $T$ isn't necessary a t-conorm. However, if $t$ is further a strictly increasing function with $t(0)=0$ then $T$ is a t-conorm. Another way is to slightly modify the function $T$ as for all $x,y\in [0,1]$,
%\begin{equation*}
%T(x,y)=\left\{
%  \begin{array}{ll}
%    \max\{x,y\} & \hbox{if }\min\{x,y\}=0, \\
%  t^{(-1)}(F(t(x),t(y))) & \hbox{otherwise.}
%  \end{array}
%\right.
%\end{equation*} Then one can check that $T$ is a t-conorm.}
%
%\item \emph{If $\infty$ be a neutral element of $F$, then $T(x,1)=t^{[-1]}(F(t(x),t(1)))\leq t^{[-1]}(t(x))$ and $t^{[-1]}(t(x))\geq x$ for all $x\in [0,1]$. So, $1$ isn't necessary a neutral element of $T$. Therefore, if $\mathfrak{T}(M)\cap (M\setminus \{t(0)\})=\emptyset$ then $T$ isn't necessary a t-norm. However, if $t$ is further a strictly increasing function with $t(1)=\infty$ then $T$ is a t-norm.}
%\end{enumerate}
%\end{remark}

\begin{example}\label{exap6.1} \emph{Let $M\in \mathcal{A}$ and $(\mathcal{U},\mathcal{V})=(\{[b_k, d_k]\mid k\in K\}, \{c_k \mid k\in \overline{K}\})$ be associated with $M$.}
\renewcommand{\labelenumi}{(\roman{enumi})}

\emph{(1) Let $F(x,y)=\max\{x,y\}$ for all $x,y\in[0,\infty]$. In Example \ref{exp4.1} (i), $M=[0,\frac{1}{4}]\cup \{\frac{1}{2}\}\cup (\frac{3}{4},1]$, and we have $\mathfrak{T_{1}}(M)=\mathfrak{T_{2}}(M)=\emptyset$ and $\mathfrak{T_{3}}(M)=(\frac{1}{2},\frac{3}{4}]$. So, the $F$-condition of $M$ holds and by Theorem \ref{th6.2}, the following function $T:[0,1]^2\rightarrow[0,1]$ given by Eq.(\ref{eq:12}) is associative:
\begin{equation*}
T(x,y)=\left\{
  \begin{array}{ll}
  \frac{3}{4}& \hbox{if }(x,y)\in[0,\frac{3}{4}]\times[\frac{1}{2},\frac{3}{4}]\cup[\frac{1}{2},\frac{3}{4}]\times[0,\frac{3}{4}], \\
   \\
   \max\{x,y\} & \hbox{otherwise.}
  \end{array}
\right.
\end{equation*}}

\emph{(2) Let $F(x,y)=x+y+xy$ for all $x,y \in [0,\infty]$ and the function $t:[0,1]\rightarrow [0,\infty]$ be defined by \begin{equation*}
 t(x)=\begin{cases}
1 & \hbox{if }\ x=0,\\
\infty  & \hbox{otherwise}.
\end{cases}
\end{equation*}
Then $M=\{1,\infty\}$, and $\mathfrak{T}(M)=\emptyset$. So, $\mathfrak{T}(M) \cap (M\setminus \{t(0)\})=\emptyset$ and by Theorem \ref{th6.2},
the following function $T:[0,1]^2\rightarrow[0,1]$ given by Eq.(\ref{eq:12}) is associative:
\begin{equation*}
T(x,y)=\left\{
  \begin{array}{ll}
  1& \hbox{if }(x,y)\in(0,1]^{2}, \\
   0 & \hbox{otherwise.}
  \end{array}
\right.
\end{equation*}}

\emph{(3) Let $F(x,y)=x+y$ for all $x,y \in [0,\infty]$ and the function $t:[0,1]\rightarrow [0,\infty]$ be defined by \begin{equation*}
 t(x)=\begin{cases}
x & \hbox{if }\ x\in[0,\frac{1}{4}],\\
\frac{1}{4} &  \hbox{if }\ x\in(\frac{1}{4},\frac{1}{2}],\\
x & \hbox{if }\ x\in(\frac{1}{2},1]
\end{cases}
\end{equation*}
Then $M=[0,\frac{1}{4}]\cup(\frac{1}{2},1]$, and we have $\mathfrak{T_{1}}(M)=\mathfrak{T_{2}}(M)=(\frac{1}{2},1]$. So, $\mathfrak{T}(M) \cap (M\setminus \{t(0)\})\neq\emptyset$ and by Theorem \ref{th6.2}, the following function $T:[0,1]^2\rightarrow[0,1]$ given by Eq.(\ref{eq:12}) isn't associative:
\begin{equation*}
T(x,y)=\left\{
  \begin{array}{ll}
   x+y & \hbox{if }0\leq x+y< \frac{1}{4} \hbox{ or }\frac{1}{2}< x+y\leq 1, \\
   \frac{1}{2} & \hbox{if }\frac{1}{4}\leq x+y\leq \frac{1}{2}, \\
   1  & \hbox{otherwise.}
  \end{array}
\right.
\end{equation*} Indeed, put $x=\frac{1}{5}, y=\frac{1}{4}, z=\frac{1}{2}$. Then $T(T(x,y),z)=T(T(\frac{1}{5},\frac{1}{4}),\frac{1}{2})=T(\frac{1}{2},\frac{1}{2})=1$, $T(x,T(y,z))=T(\frac{1}{5},T(\frac{1}{4},\frac{1}{2}))=T(\frac{1}{5},\frac{3}{4})=\frac{19}{20}$. Thus, $T$ isn't an associative function.}
\end{example}

\section{Further properties of the function $T$ given by \eqref{eq:12}}
In this section, we explore the idempotence, the limit property, the conditional cancellation law and the continuity of the function $T$ given by \eqref{eq:12}, respectively.

Denote by $\varGamma$ the set of all functions $F:[0,\infty]^2\rightarrow[0,\infty]$ such that $([0,\infty],F,\leq)$ is a fully ordered Abel semigroup with $F(x,0) \geq x$ for all $x \in [0,\infty]$.

In what follows, we always suppose that $t:[0,1]\rightarrow [0,\infty]$ is a left continuous non-decreasing function and  $F\in \varGamma$.

\subsection{Idempotence}
Recall that a function $T:[0,1]^2\rightarrow[0,1]$ is said to be idempotent if $T(x,x)=x$ for all $x\in[0,1]$. In this subsection we investigate the idempotence of $T$ given by \eqref{eq:12}.

%Denote
%$$\mathbb{H}=\{c\mid \mbox{ there are an } x_0\in[0,1]\mbox{ and }\varepsilon >0\mbox{ such that }t|_{[x_0,x_0+\varepsilon]}=c\},$$
%$$\mathbb{N}=\{\sup \{x\in[0,1]\mid t(x)=y\}\mid y\in \mathbb{H}\},~\mathbb{W}=\{x\in [0,1]\mid t(x)\in M\setminus \mathbb{H} \},$$
%$$\mathbb{D}=\mathbb{N}\cup \mathbb{W}.$$
%
%In particular, $t^{[-1]}(x) \in \mathbb{D}$ for all $x\in [0,\infty]$. Moreover, $t^{[-1]}\circ t(x)=x$ for all $x\in \mathbb{D}$.
%and  $t \circ t^{[-1]}(x)=x$ for all $x\in M$. In  particular, we write $\mathbb{R}=\{x\in [0,1]\mid t(x)\in M\setminus \mathbb{H} \mbox{ and } t(x^{+})\notin \mathbb{H}\}.$

We first have the proposition as below.
\begin{proposition}\label{prop7.1}
Let the function $T$ be given by Eq.\eqref{eq:12}, $t:[0,1]\rightarrow [0,\infty]$ be a left continuous non-decreasing function with $M\in \mbox{Ran}(t)$ and $F\in \varGamma$.
Then for any $x\in[0,1)$, $T(x,x)=x$ if and only if  $x \in \mathbb{D}$ and $M\cap [t(x),F(t(x),t(x))]=\{t(x)\}$.
\end{proposition}
\begin{proof}Let $x\in[0,1)$ be such that $T(x,x)=x$. Assume $x \notin \mathbb{D}$. Then there exists an $x_{1}\in \mathbb{G}$ such that $x<x_{1}$, $t(x_{1})=t(x)$ and $t^{[-1]}(t(x))=x_{1}$. Thus $T(x,x)=t^{[-1]}(F(t(x),t(x)))\geq t^{[-1]}(t(x))> x$, a contradiction.  Therefore,  $x \in \mathbb{D}$.  Moreover, we have
$t^{[-1]}(F(t(x),t(x)))= T(x,x)=x=t^{[-1]}(t(x))$. Since $F(t(x),t(x))\geq t(x)$, we have either $F(t(x),t(x))= t(x)$ when $t(x)=t(x^+)$ or $F(t(x),t(x))\leq t(x^+)$ when $t(x)\neq t(x^+) \notin M$ or $F(t(x),t(x))< t(x^+)$ when $t(x)\neq t(x^+) \in M$.
Thus, $M\cap [t(x),F(t(x),t(x))]=\{t(x)\}$.

The converse is obvious.
\end{proof}

Further, we have the following corollary.
\begin{corollary}\label{cor7.01}
	Let the function $T$ be given by Eq.\eqref{eq:12}. If $t:[0,1]\rightarrow [0,\infty]$ is a left continuous strictly increasing function,
	then for any $x\in[0,1)$, $T(x,x)=x$ if and only if $M\cap [t(x),F(t(x),t(x))]=\{t(x)\}$.
\end{corollary}

 A function $F:[0,\infty]^2\rightarrow[0,\infty]$ is said to be strictly increasing if $F(x,y)<F(x,z)$ whenever $\infty >x$ and $y<z$. Then from Proposition \ref{prop7.1} we have the following remark.
\begin{remark}\label{remk7.1}
	\emph{Let the function $T$ be given by Eq.~\eqref{eq:12}, $t:[0,1]\rightarrow [0,\infty]$ be a left continuous non-decreasing function with $M\in \mbox{Ran}(t)$ and $F\in \varGamma$.
		\renewcommand{\labelenumi}{(\roman{enumi})}
		\begin{enumerate}
			\item Obviously, $1$ is an idempotent element of $T$.
		\item If there exists an $x\in\mathbb{D}$ such that $t(x)$ is an idempotent element of $F$, then $x$ is an idempotent element of $T$.
		\item If there exists an $x\in\mathbb{D}$ such that $t(x^+)$ is an idempotent element of a strictly increasing function $F$, then $x$ is an idempotent element of $T$.
        \item If $t$ is continuous at $x\in\mathbb{D}$, then $T(x,x)=x$ if and only if $F(t(x),t(x))=t(x)$.
		\end{enumerate}}
\end{remark}

\begin{example}\label{exp04.1}\emph{Let $t:[0,1]\rightarrow [0,\infty]$ be defined by
		    		\begin{equation*}
    			t(x)=\begin{cases}
    				5x & \ \hbox{if } x\in[0,0.2],\\
    				2 &  \ \hbox{if } x\in(0.2,0.5],\\
    				10x &  \ \hbox{if } x\in(0.5,1]\\
    			\end{cases}
    		\end{equation*}
		and the function $F\colon[0,\infty]^{2}\to[0,\infty]$ be defined by $F(x,y)=x+ y$. Then by Eq.~\eqref{eq:12}, \begin{equation*}
    			     T(x,y)=\begin{cases}
    				x+y & \ \hbox{if }x+y\leq0.2,\\
    				0.2 & \ \hbox{if }(x,y)\in [0,0.2)^2 \setminus\{x+y\leq 0.2\},\\
    				0.5 & \ \hbox{if }(x,y)\in [0,0.5]^2 \setminus [0,0.2)^2,\\
    				0.5x+y &\ \hbox{if } (x,y)\in ([0.0.2]\times[0.5,1])\cap\{ 1<x+2y<2\},\\
    			y+0.2	& \ \hbox{if } (x,y)\in (0.2,0.5]\times[0.5,0.8),\\
    			0.5y+x & \ \hbox{if } (x,y)\in ([0.5,1]\times[0.0.2])\cap\{ 1<2x+y<2\},\\
    			x+0.2	& \ \hbox{if } (x,y)\in [0.5,0.8)\times(0.2,0.5],\\
    			1 & \ \mbox{otherwise}.
    			\end{cases}
    		\end{equation*}
		By Proposition \ref{prop7.1}, one can check that 0, 0.5,1 are idempotent elements of the function $T$, respectively.}
\end{example}

\begin{theorem}
Let $t:[0,1]\rightarrow [0,\infty]$ be a left continuous non-decreasing function with $M\in \mbox{Ran}(t)$ and $F\in \varGamma$. Then the function $T$ given by Eq.~\eqref{eq:12} is idempotent
 if and only if $t$ is strictly increasing on $[0, 1]$ and $F(t(x),t(x))\leq t(x^+)$ for all $x\in[0,1)$.
\end{theorem}
\begin{proof}If the function $T$ given by Eq.\eqref{eq:12} is idempotent, i.e., $T(x, x) = x$ for all $x\in[0,1]$, then from Proposition \ref{prop7.1}, $x\in \mathbb{D}$, thus $\mathbb{D}=[0, 1)$. Therefore, $t$ is strictly increasing on $[0, 1]$. The rest of the proof is completely analogous to Proposition 3.1 of \cite{Wang}.
\end{proof}

\subsection{Limit property}
The function $T$ given by Eq.\eqref{eq:12} is said to have the limit property if $\lim_ {n \rightarrow \infty} x_{T}^{(n)}=1$ for all $x\in (0,1)$. In this subsection,
 we investigate some necessary and sufficient conditions for the function $T$ given by Eq.\eqref{eq:12} satisfying the limit property.

\begin{proposition}\label{prop6.3}
Let $t:[0,1]\rightarrow [0,\infty]$ be a left continuous non-decreasing function with $M\in \mbox{Ran}(t)$ and $F\in \varGamma$. If the function $T$ given by Eq.\eqref{eq:12} has the limit property, then $F(t(x),t(x))\geq t(x^+)$ for all $x\in(0,1)$.
\end{proposition}
\begin{proof}Assuming there exists an $x\in(0,1)$ such that $F(t(x),t(x))< t(x^+)$.  We write $k=t^{[-1]}(t(x))$. Then $x\leq k<1$ and $t(k)=t(x)$, so that $k= t^{[-1]}(F(t(x),t(x)))=x_{T}^{(2)}$. Suppose that $k=x_{T}^{(n)}$ when $n=m$.  Then $F(t(x),t(k))= F(t(x),t(x_{T}^{(m)}))$ when $n=m+1$. Moreover, $t(x^+)> F(t(x),t(k))= F(t(x),t(x_{T}^{(m)}))$. This means $k= t^{[-1]}(F(t(x),t(x_{T}^{(m)})))=x_{T}^{(m+1)}$. Hence, by a mathematical induction, $k=x_{T}^{(n)}$ for all $n\in\{2,3,\cdots\}$. Therefore, $1>k= \lim_ {n \rightarrow \infty} x_{T}^{(n)}$, a contradiction.
\end{proof}

\begin{proposition}\label{prop6.4}
Let $t:[0,1]\rightarrow [0,\infty]$ be a left continuous non-decreasing function with $M\in \mbox{Ran}(t)$ and $F\in \varGamma$ be continuous. If the function $T$ given by Eq.\eqref{eq:12} has the limit property, then $F(t(0^{+}),t(x))\geq t(x^+)$ for all $x\in(0,1)$.
\end{proposition}
\begin{proof}Assuming there exists a $y\in(0,1)$ such that $F(t(0^{+}),t(y))< t(y^+)$. It follows that there exists an $x\in(0,y)$ near enough to 0 such that $F(t(x),t(y))< t(y^+)$ by the continuity of $F$. We write $k=t^{[-1]}(t(y))$. Then $y\leq k<1$. Since  $t(y)\leq F(t(x),t(y))< t(y^+)$, we have  $k= t^{[-1]}(F(t(x),t(y)))\geq t^{[-1]}(F(t(x),t(x)))=x_{T}^{(2)}$. Suppose that $k\geq x_{T}^{(n)}$ when $n=m$. Then $F(t(x),t(k))\geq F(t(x),t(x_{T}^{(m)}))$  when $n=m+1$. Since $t(k)=t(y)$, we have $t(k^+)> F(t(x),t(k))\geq F(t(x),t(x_{T}^{(m)}))$. This means $k= t^{[-1]}(F(t(x),t(k)))\geq t^{[-1]}(F(t(x),t(x_{T}^{(m)})))=x_{T}^{(m+1)}$. Hence, by a mathematical induction, $k\geq x_{T}^{(n)}$ for all $n\in\{2,3,\cdots\}$. Therefore, $1>k\geq \lim_ {n \rightarrow \infty} x_{T}^{(n)}$, a contradiction.
\end{proof}

\begin{proposition}\label{prop6.5}
Let $t:[0,1]\rightarrow [0,\infty]$ be a left continuous non-decreasing function with $M\in \mbox{Ran}(t)$ and $F\in \varGamma$ be continuous. If  $F(t(0^{+}),t(x))> t(x^+)$ for all $x\in(0,1)$, then the function $T$ given by Eq.\eqref{eq:12} has the limit property.
\end{proposition}
\begin{proof}Suppose that there exists an $x\in(0, 1)$ such that $\lim_{n\rightarrow\infty} x^{n}_{T}=y< 1$, i.e., the non-decreasing sequence $\{x^{n}_{T}\}_{n\in N}$ converges to $y$. Then $y\geq x^{n}_{T}$ for all $n\in N$. First, we prove that there is an $m \in N$ such that $x^{n}_{T}=y$ for all $n > m$. Indeed, $t(y))\leq F(t(x),t(y))$ since $F\in \varGamma$. This implies that there exists an $\epsilon>0$ such that $F(t(x),t(z))\geq t(y)$ for all $z\in(y-\epsilon,y]$ since $F$ is continuous. Thus there exists an $m \in N$ such that $x^{n}_{T} \in (y-\epsilon,y]$ for all $n \geq m$, which means $F(t(x),t(x_{T}^{(n)}))\geq t(y)$ for all $n>m$. Hence, from the monotonicity of $t^{[-1]}$ we have $y\geq x_{T}^{(n+1)}=T(x,x_{T}^{(n)})=t^{[-1]}(F(t(x),t(x_{T}^{(n)})))\geq t^{[-1]}(t(y))\geq y$. Therefore, there is an $m\in N$ such that $x_{T}^{(n)}=y$ for all $n>m$. This follows that $y=x_{T}^{(n+1)}=T(x,x_{T}^{(n)})=T(x,y)=t^{[-1]}(F(t(x),t(y)))$ for all $n>m$, which deduces that $F(t(x),t(y))\leq t(y^+)$, contrary to the fact $F(t(x),t(y))\geq F(t(0^{+}),t(y))> t(y^+)$.
\end{proof}

\begin{theorem}\label{Thm6.6}
Let $t:[0,1]\rightarrow [0,\infty]$ be a left continuous non-decreasing function and $t^{[-1]}(t(x))=t^{[-1]}(t(x^+))$ for all $x\in(0,1)$.
If $F\in \varGamma$ is continuous, then the following assertions are equivalent:
\renewcommand{\labelenumi}{(\roman{enumi})}
\begin{enumerate}
\item $F(t(0^{+}),t(x))> t(x^+)$ for all $x\in(0,1)$.
\item The function $T$ given by Eq.\eqref{eq:12} has the limit property.
\end{enumerate}
\end{theorem}
\begin{proof}From Proposition \ref{prop6.5}, (i) implies (ii). In order to prove that (ii) implies (i), assume there exists a $y\in(0,1)$ such that $F(t(0^{+}),t(y))\leq t(y^+)$. It follows that there exists an $x\in(0,y)$ near enough to 0 such that $F(t(x),t(y))\leq t(y^+)$ by the continuity of $F$. Put $k=t^{[-1]}(t(y))=t^{[-1]}(t(y^+))$. Then $y\leq k<1$ and $t(k)=t(y)$.
 Because of $t(y)\leq F(t(x),t(y))\leq t(y^+)$, we have  $k\geq t^{[-1]}(F(t(x),t(y)))\geq t^{[-1]}(F(t(x),t(x)))=x_{T}^{(2)}$. Suppose that $y\geq x_{T}^{(n)}$  when $n=m$. Then $F(t(x),t(y))\geq F(t(x),t(x_{T}^{(m)}))$  when $n=m+1$. We have $t(y^+)\geq F(t(x),t(y))\geq F(t(x),t(x_{T}^{(m)}))$ since $t(k)=t(y)$. Because of $k=t^{[-1]}(t(y))=t^{[-1]}(t(y^+))$, we have $k=t^{[-1]}(F(t(x),t(k))) \geq t^{[-1]}(F(t(x),t(x_{T}^{(m)})))=x_{T}^{(m+1)}$. Hence, by a mathematical induction, $k\geq x_{T}^{(n)}$ for all $n\in\{2,3,\cdots\}$. Therefore, $1>k\geq \lim_ {n \rightarrow \infty} x_{T}^{(n)}$, a contradiction.
\end{proof}

Notice that the condition $t^{[-1]}(t(x))=t^{[-1]}(t(x^+))$ in Theorem \ref{Thm6.6} cannot be deleted generally.

\begin{example}\label{exp4.01}\emph{Let $t:[0,1]\rightarrow [0,\infty]$ be defined by
		    		\begin{equation*}
    			t(x)=\begin{cases}
    				1 & \ \hbox{if } x\in[0,0.5],\\
    				2 &  \ \hbox{if } x\in(0.5,1]\\
    			\end{cases}
    		\end{equation*}
		and the function $F\colon[0,\infty]^{2}\to[0,\infty]$ be defined by $F(x,y)=x+y$. Then by Eq.~\eqref{eq:12}, \begin{equation*}
    			     T(x,y)=1 \mbox{ for all } (x,y)\in[0,1]^2.
    		\end{equation*}
		It is easy to check that $t:[0,1]\rightarrow [0,\infty]$ is a left continuous non-decreasing function, $F\in \varGamma$ is continuous and $T$ has the limit property. However, $F(t(0^{+}),t(0.5))=2=t(0.5^+)$.}
\end{example}

From Theorem \ref{Thm6.6}, we have the following corollary.
\begin{corollary}
Let $t:[0,1]\rightarrow [0,\infty]$ be a  left continuous strictly increasing function. If $F\in \varGamma$ is continuous, then the following assertions are equivalent:
\renewcommand{\labelenumi}{(\roman{enumi})}
\begin{enumerate}
\item $F(t(0^{+}),t(x))> t(x^+)$ for all $x\in(0,1)$.
\item The function $T$ given by Eq.\eqref{eq:12} has the limit property.
\end{enumerate}
\end{corollary}

\begin{proposition}\label{prop6.7}
Let $t:[0,1]\rightarrow [0,\infty]$ be a left continuous non-decreasing function and $t^{[-1]}(t(x))=t^{[-1]}(t(x^+))$ for all $x\in(0,1)$.
If $F\in \varGamma$ is strict, then the following assertions are equivalent:
\renewcommand{\labelenumi}{(\roman{enumi})}
\begin{enumerate}
\item $F(t(0^{+}),t(x))\geq t(x^+)$ for all $x\in(0,1)$.
\item The function $T$ given by Eq.\eqref{eq:12} has the limit property.
\end{enumerate}
\end{proposition}
\begin{proof}From  Proposition \ref{prop6.4}, (ii) implies (i).
In order to prove that (i) implies that (ii), assume  there exists an $x\in(0,1)$ such that $\lim_ {n \rightarrow \infty} x_{T}^{(n)}=y<1$. Then by the proof of Proposition \ref{prop6.5} we have $F(t(x),t(y))\leq t(y^+)$. On the other hand, $F(t(x),t(y))> F(t(0^{+}),t(y))\geq t(y^+)$ since $F\in \varGamma$ is strict, a contradiction.
\end{proof}

Generally, the condition $t^{[-1]}(t(x))=t^{[-1]}(t(x^+))$ in Proposition \ref{prop6.7} cannot be removed.
\begin{example}\label{exp4.001}\emph{Let $t:[0,1]\rightarrow [0,\infty]$ be defined by
		    		\begin{equation*}
    			t(x)=\begin{cases}
    				1 & \ \hbox{if } x\in[0,0.5],\\
    				4x &  \ \hbox{if } x\in(0.5,1]\\
    			\end{cases}
    		\end{equation*}
		and the function $F\colon[0,\infty]^{2}\to[0,\infty]$ be defined by $F(x,y)=x+y$. Then by Eq.~\eqref{eq:12}, \begin{equation*}
    			     T(x,y)=\begin{cases}
    				0.5 &  \ \hbox{if } (x,y)\in [0,0.5]^2,\\
    				0.25+y & \ \hbox{if } (x,y)\in [0.0.5]\times(0.5,0.75],\\
    			0.25+x & \ \hbox{if } (x,y)\in (0.5,0.75]\times[0.0.5],\\
    			1 & \ \mbox{otherwise}.
    			\end{cases}
    		\end{equation*}
		It is easy to check that $t:[0,1]\rightarrow [0,\infty]$ is a left continuous non-decreasing function, $F\in \varGamma$ is strict and $F(t(0^{+}),t(x))\geq t(x^+)$ for all $x\in(0,1)$. However, $\lim_ {n \rightarrow \infty} 0.5_{T}^{(n)}=0.5$, i.e., $T$ has not the limit property.}
\end{example}

From Proposition \ref{prop6.7}, we have the following corollary.
\begin{corollary}
Let $t:[0,1]\rightarrow [0,\infty]$ be a  left continuous strictly increasing function. If $F\in \varGamma$ is strict, then the following assertions are equivalent:
\renewcommand{\labelenumi}{(\roman{enumi})}
\begin{enumerate}
\item $F(t(0^{+}),t(x))\geq t(x^+)$ for all $x\in(0,1)$.
\item The function $T$ given by Eq.\eqref{eq:12} has the limit property.
\end{enumerate}
\end{corollary}

In particular, we have the following corollary.
\begin{corollary}
Let  $t:[0,1]\rightarrow [0,\infty]$ be a left continuous strictly increasing function. If $F\in \varGamma$ is strict and $t(0^{+})$ is a neutral element of $F$, then the following assertions are equivalent:
\renewcommand{\labelenumi}{(\roman{enumi})}
\begin{enumerate}
\item $t$ is continuous.
\item The function $T$ given by Eq.\eqref{eq:12} has the limit property.
\end{enumerate}
\end{corollary}

\begin{proposition}\label{prop6.8}
Let $t:[0,1]\rightarrow [0,\infty]$ be a left continuous non-decreasing function and $F\in \varGamma$ be continuous. If the function $T$ given by Eq. \eqref{eq:12} is a triangular supconorm, then the following assertions are equivalent:
\renewcommand{\labelenumi}{(\roman{enumi})}
\begin{enumerate}
\item $T(x, x) > x$ for all $x \in(0, 1)$.
\item The function $T$ given by Eq.\eqref{eq:12} has the limit property.
\end{enumerate}
\end{proposition}
\begin{proof}Obviously, (ii) implies (i). In order to show that (i) $\Rightarrow$ (ii), suppose that there exists an $x\in(0, 1)$ such that $\lim_{n\rightarrow\infty} x^{n}_{T}=y< 1$.
Similar to the proof of Proposition \ref{prop6.5} we can prove that there exits an $m \in N$ such that $ x^{n}_{T}= y$ for all $n > m$.
 Then $T(y, y) = T(x^{m+1}_{T}, x^{m+1}_{T}) = x^{2m+2}_{T}= y$ since $T$ is associative. Therefore, $y = 0$ by (i), contrary to $y \in(0, 1)$.
\end{proof}

\begin{proposition}
 Let $t:[0,1]\rightarrow [0,\infty]$ be a left continuous non-decreasing function and $F\in \varGamma$ be continuous. If the function $T$ given by Eq. \eqref{eq:12} is a triangular supconorm which satisfies the conditional cancellation law, then
the function $T$ given by Eq.\eqref{eq:12} has the limit property.
\end{proposition}
\begin{proof}Suppose that $\lim_{n\rightarrow\infty} x^{n}_{T}< 1$ for an $x \in(0, 1)$. Then by Proposition \ref{prop6.8} there is a $y \in(0, 1)$ such that $T (y, y) = y$. Thus  $y\leq T(y,0)\leq T(y, y)=y$, i.e., $T(y,0)=T(y, y)=y$, then $y=0$ since $T$ satisfies the conditional cancellation law, contrary to $y \in(0, 1)$.
\end{proof}

\subsection{The conditional cancellation law}
The function $T$
 defined by Eq. \eqref{eq:12} is said to satisfy the conditional cancellation law if $T(x,y)=T(x,z)<1$ for all $x,y,z\in[0,1]$ implies $y=z$. The function $T$ is said to satisfy the cancellation law if $T(x,y)=T(x,z)$) for all $x,y,z\in[0,1]$ implies  $x<1$ or $y=z$.
In this subsection, we obtain a necessary and sufficient condition for the function $T$
defined by Eq. \eqref{eq:12} satisfying the conditional cancellation law, which is applied for showing necessary and sufficient conditions for the function $T$ to be a t-supconorm.

%\begin{example}\label{ex3.4}
%\emph{Let $t:[0,1]\rightarrow [0,\infty]$ be defined by \begin{equation*}
%  		t(x)=\begin{cases}
%  			2 & \ \hbox{if } x\in[0,0.5],\\
%  			4x &  \ \hbox{if } x\in(0.5,1]
%  		\end{cases}
%  	\end{equation*} and the function $F\colon[0,\infty]^{2}\to[0,\infty]$ be defined by $F(x,y)=x+y+xy$. Then by Eq.\eqref{eq:12},
%  		$$T(x,y)=1 \mbox{ ~for all } (x,y)\in[0,1]^2. $$
%  	Obviously, $T$ is a t-superconorm which satisfies the conditional cancellation law.}
%\end{example}
%
% Example \ref{ex3.4} show that the function $T$ given by Eq. \eqref{eq:12} satisfies the conditional cancellation law, then $t\colon[0,1]\to[0,\infty]$ does not have to be strictly increasing.

First note that if $t(0)\in \mathbb{H}$ then one can check that the function $T$ defined by Eq. \eqref{eq:12} satisfies the conditional cancellation law if and only if $T(x,y)=1$ for all
$(x,y)\in [0,1]^2$. Therefore, in this subsection, we always assume that $t\colon[0,1]\to[0,\infty]$ is a left continuous non-decreasing function with $t(0)\notin \mathbb{H}$ and $F\in \varGamma$ is strict. Let $\beta=\min \mathbb{H}$. Then there exists an $\alpha \in [0,1]$ such that $t(\alpha)=\beta$ and $t(\alpha-\varepsilon)< \beta$ for an arbitrary $\varepsilon \in (0,1-\alpha)$. Moreover, $t\mid_{[0, \alpha]}$ is strictly increasing.

We have the following lemma.
\begin{lemma}\label{lem3.00} Let $(\mathcal{U},\mathcal{V})=(\{[b_k, d_k]\mid k\in K\}, \{c_k \mid k\in \overline{K}\})$ be associated with $M\in \mathcal{A}$ and $F\in \varGamma$ be strict. If $F(\mathbb{H},M)\in [t(1),\infty]$, then the following are equivalent:
\renewcommand{\labelenumi}{(\roman{enumi})}
\begin{enumerate}
\item $F((M\setminus C), M)\subseteq (M\setminus C) \cup [t(1),\infty]$
\item $F((M\setminus C), M)\subseteq M \cup [t(1),\infty]$.
\end{enumerate}
\end{lemma}
\begin{proof}It is trivial that (i) implies (ii). In order to show that (ii) implies (i), suppose that $F((M\setminus C),M)\nsubseteq (M\setminus C) \cup [t(1),\infty]$.
 Then there exist an $x \in M\setminus C, y \in M$ and $k \in K$ such that $F(x,y) \in [b_{k}, d_{k}]$. On the other hand, from $F(x,y)\in M\cup [t(1),\infty]$ we have $F(x,y)= b_{k}$ where $\{b_{k}\}=M \cap [b_{k}, d_{k}]$. Since $x \in (M\setminus C) \cap [0,t(1))$,
 there exists an $\varepsilon >0$ such that $a \in M\setminus C$ for an arbitrarily $a\in[x,x+\varepsilon]$. Furthermore, by (ii), there exists an $l\in K$ such that $F(a,y)= b_{l}$ where $\{b_{l}\}=M\cap[b_{l}, d_{l}]$.
 Clearly, $b_{l}$ has different values when $a$ changes since $F$ is strict.  Therefore, $\{[b_{l}, d_{l}]|l\in K\mbox{ and }b_{l}\in [b_{l}, d_{l}]\}$ is an uncountable set since $[x,x+\varepsilon]$ is uncountable, contrary to the fact that $K$ is a countable set.
\end{proof}

\begin{theorem}\label{Thm4.3}
Let $(\mathcal{U},\mathcal{V})=(\{[b_k, d_k]\mid k\in K\}, \{c_k \mid k\in \overline{K}\})$ be associated with $M\in \mathcal{A}$ and $F\in \varGamma$ be strict. Then the function $T$ given by Eq.~\eqref{eq:12} is conditionally cancellative if and only if the following  hold:
\renewcommand{\labelenumi}{(\roman{enumi})}
\begin{enumerate}
\item $F(M\setminus C, M)\subseteq  M \cup [t(1),\infty]$.
\item $F(\mathbb{H},M)\in [t(1),\infty]$.
\end{enumerate}
\end{theorem}
\begin{proof}We prove that the function $T$ given by Eq.~\eqref{eq:12} is not
conditionally cancellative if and only if either $F(M\setminus C, M)\nsubseteq  M \cup [t(1),\infty]$ or $F(\mathbb{H},M)\notin [t(1),\infty]$.

Suppose that $T$ is not conditionally cancellative, i.e., there exist three elements $x_{1}, x_{2}, y \in [0,1]$ with $x_{1} < x_{2}$ and $t(y)<1$ such that
$T(x_{1}, y) =T(x_{2},y) < 1$. In the following we distinguish two cases.

(1) If $F(\mathbb{H},M)\in [t(1),\infty]$, then $t(x_{1}), t(x_{2})\notin \mathbb{H}$ since $T(x_{1}, y) =T(x_{2},y) < 1$. Thus, $F(t(x_{1}),t(y))<F(t(x_{2}),t(y))$ since $F$ is strict. It follows from $T(x_{1}, y) =T(x_{2},y)$ that \begin{equation}\label{eq001}(M\setminus C)\cap[F(t(x_{1}),t(y)),F(t(x_{2}),t(y))]=\emptyset.\end{equation}  On the other hand, $(M\setminus C)\cap(t(x_{1}), t(x_{2}))\neq \emptyset$ since $t(x_{1}), t(x_{2})\notin \mathbb{H}$. This follows that there is an $x \in [0,1]$ such that $$t(x) \in (M\setminus C)\cap(t(x_{1}), t(x_{2})).$$ Clearly, $F(t(x),t(y))\in [F(t(x_{1}),t(y)),F(t(x_{2}),t(y))]$ since $F$ is strict, which together with \eqref{eq001} implies $F(t(x),t(y))\notin (M\setminus C)\cup [t(1),\infty]$. Therefore, from Lemma \ref{lem3.00} we have $F(M\setminus C, M)\nsubseteq M\cup [t(1),\infty]$.

(2) If $F(M\setminus C, M)\subseteq  M \cup [t(1),\infty]$ then, by (1), we immediately have that $F(\mathbb{H},M)\notin [t(1),\infty]$.

Conversely, if $F(M\setminus C, M)\nsubseteq  M \cup [t(1),\infty]$, i.e., there exist three elements $x,y \in [0,1]$, $k\in K$ such that $t(x)\in M\setminus C$, $t(y) \in M$ and $F(t(x),t(y)) \in [b_{k}, d_{k}]\setminus\{c_{k}\}$ where $d_{k}<t(1)$ and $\{c_{k}\}=M\cap[b_{k}, d_{k}]$, then there exist two elements $t(x_{1}), t(x_{2}) \in M$ with $t(x_{1}) \neq t(x_{2})$ such
that $F(t(x_{1}), t(y)), F(t(x_{2}), t(y)) \in [b_{k}, d_{k}]\setminus\{c_{k}\}$ since $t(x) \in (M\setminus C) \cap [0,\infty]$. Therefore, $T(x_{1}, y)=T(x_{2}, y)<1$, i.e., $T$ is not conditionally cancellative.

If $F(\mathbb{H},M)\notin [t(1),\infty]$, then there exist two elements $x_{1}, x_{2}\in [0,1]$ with $t(x_{1})=t(x_{2})$ such that $F(t(x_{1}),t(y))= F(t(x_{2}),t(y))<t(1)$ for all $y \in [0,1]$. Therefore, $T(x_{1}, y)=T(x_{2}, y)<1$, i.e., $T$ is not conditionally cancellative.
\end{proof}

\begin{corollary}\label{cor4.3}
Let $(\mathcal{U},\mathcal{V})=(\{[b_k, d_k]\mid k\in K\}, \{c_k \mid k\in \overline{K}\})$ be associated with $M\in \mathcal{A}$ and $F\in \varGamma$ be strict. Then the function $T$ given by Eq.~\eqref{eq:12} is cancellative if and only if the following hold:
\renewcommand{\labelenumi}{(\roman{enumi})}
\begin{enumerate}
\item $t:[0,1]\to[0,\infty]$ is a left continuous strictly increasing function.
\item $F(M\setminus C, M)\subseteq M$.
\end{enumerate}
\end{corollary}
\begin{proof}Let the function $T$ given by Eq.~\eqref{eq:12} be cancellative. Then, obviously, $T$ is conditionally cancellative and $t:[0,1]\to[0,\infty]$ is strictly increasing with $t(1)=\infty$. Thus $[t(1),\infty]\subseteq  M$ and $F(\mathbb{H},M)=\emptyset$. Therefore, from Theorem \ref{Thm4.3}, we have $F(M\setminus C, M)\subseteq  M \cup [t(1),\infty]$ and $F(\mathbb{H},M)\in [t(1),\infty]$, which means $F(M\setminus C, M)\subseteq  M$.

Conversely, let $F(M\setminus C, M)\subseteq  M$. Then $F(x, t(1))=t(1)$ for all $x\in[0,\infty]$, thus $t(1)=\infty$. We get $T(M\setminus C, M)\subseteq  M \cup [t(1),\infty]$ and $F(\mathbb{H},M)\in [t(1),\infty]$ since $t$ is strictly increasing with $t(1)=\infty$. Therefore, from Theorem \ref{Thm4.3} $T$ is conditionally cancellative, moreover, $T$ is cancellative since $t(1)=\infty$.
\end{proof}

Lemma \ref{lem6.1} and Proposition \ref{prop5.3} imply the following lemma.
\begin{lemma}\label{lem3.5}
Let $t:[0,1]\rightarrow [0,\infty]$ be a left continuous non-decreasing function with $\mbox{Ran}(t)=M$. Then for any $x,y\in[0,\infty]$, $t^{[-1]}(x)=t^{[-1]}(y)$ if and only if $(\min\{x,y\},\max\{x,y\}]\cap (M\setminus t(0))=\emptyset$.
\end{lemma}
%\begin{proof}Let for any $x,y\in[0,\infty]$, $t^{[-1]}(x)=t^{[-1]}(y)$. Then $t\circ t^{[-1]}(x)= t\circ t^{[-1]}(y)$. From Lemma \ref{lem6.1} and Proposition \ref{prop5.3}, it is immediate to $(\min\{x,y\},\max\{x,y\}]\cap (M\setminus t(0))=\emptyset$.
%
%Conversely, If $(\min\{x,y\},\max\{x,y\}]\cap (M\setminus t(0))=\emptyset$, without loss of generality, we say $x< y$, then $(x,y]\subseteq (0,t(0)]$ or there is $k\in K$ such that $(x,y]\subseteq [b_k,d_k]\setminus\{c_{k}\}$. If $(x,y]\subseteq (0,t(0)]$ then $t^{[-1]}(x)=0=t^{[-1]}(y)$. If $(x,y]\subseteq [b_k,d_k]\setminus\{c_{k}\}$ then $t^{[-1]}(x)=t^{[-1]}(b_k)=t^{[-1]}(y)$.
%\end{proof}

Let $M\subseteq [0,\infty]$, $K_{1}=K\cup \{\kappa\}$. Define $$H_{\kappa}=O(\{t(0)\}\cup \{z\in[0,t(0))\mid  \mbox{ there exit two elements } x,y\in M \mbox{ such that } F(x,y)=z\})$$ and  $$H_{k}=O(\{b_{k}\}\cup \{z\in(b_{k},d_{k})\mid  \mbox{ there exit two elements } x,y\in M \mbox{ such that }  F(x,y)=z\})$$for each $k \in K$.  Then we have the following theorem.

\begin{theorem}Let $(\mathcal{U},\mathcal{V})=(\{[b_k, d_k]\mid k\in K\}, \{c_k \mid k\in \overline{K}\})$ be associated with $M\in \mathcal{A}$ and $F\in \varGamma$ which is strict with $t(0)$ a neutral element. If the function $T$ given by Eq.~\eqref{eq:12} is conditionally cancellative, then the following are equivalent:
\renewcommand{\labelenumi}{(\roman{enumi})}
\begin{enumerate}
\item $T$ is a t-supconorm.
\item $F(\cup_{k\in K_{1}}H_{k}, M)\cap (M\setminus t(0))= \emptyset$.
\end{enumerate}
\end{theorem}
\begin{proof}(i)$\Rightarrow$(ii) If $F(M , M) \subseteq M \cup [t(1),\infty]$, then either $\cup_{\alpha\in K_{1}}H_{k}=\emptyset$ or $\cup_{k\in K_{1}}H_{k}\subseteq(t(1),\infty]$. Thus $F(\cup_{k\in K_{1}}H_{k}, M)\cap (M\setminus t(0))= \emptyset$.

Now, let $F(M , M) \nsubseteq M \cup [t(1),\infty]$. Then there exist two elements $x, y \in [0,1]$ such that $F(t(x),t(y)) \notin M \cup [t(1),\infty]$, which means that $t(0)<\min\{t(x),t(y)\}<\beta=\min \mathbb{H}$ and there is a $k \in K$ such that $F(t(x),t(y))\in [b_{k}, d_{k}]\setminus\{c_{k}\}$ where $\{c_{k}\}=[b_{k}, d_{k}] \cap M$ and $d_{k}<t(1)$. On the other hand, $(M\setminus C)\setminus \{\infty\}\neq \emptyset$ since $t\mid_{[0, \alpha]}$ is strictly increasing, which implies that one can choose a $t(z) \in (M\setminus C)\setminus \{\infty\}$ such that $t(z) < \min\{t(x),t(y)\}$. Thus $F(t(y),t(z))<F(t(y),t(x))<t(1)$ since $F$ is strict. Meanwhile, from Theorem \ref{Thm4.3}(i) and Lemma \ref{lem3.00} we have $$F(M\setminus C, M) \subseteq (M\setminus C) \cup [t(1),\infty] \mbox{ and } F(\mathbb{H},M)\subseteq[t(1),\infty].$$
Hence $F(t(y),t(z))\in M\setminus C$. Therefore, by Eq.~\eqref{eq:12} we have
$T(x,T(y,z))=t^{[-1]}(F(t(x), F(t(y),t(z))))$ and $T(T(x,y),z)=t^{[-1]}(F(b_{k}, t(z)))$.
Consequently, by the associativity of both $T$ and $F$ we have
$$t^{[-1]}(F(F(t(x),t(y)),t(z)))=t^{[-1]}(F(b_{k}, t(z))),$$
implying $M\cap[F(b_{k}, t(z)), F(F(t(x),t(y)),t(z))]$ contains at most one element. In the following we prove that $t(1) \leq \min\{F(b_{k}, t(z)), F(F(t(x),t(y)),t(z))\}<\infty$.

Indeed, if $F(b_{k}, t(z))< t(1)$ then $F(b_{k}, t(z))\in M\setminus C$, implying $F(F(t(x),t(y)),t(z))=F(b_{k}, t(z))$. Since $t(z) < \infty$ and $F\in \varGamma$ is strict, we have $F(t(x), t(y))=b_{k} $, contrary to $F(t(x),t(y))\in [b_{k}, d_{k}]\setminus\{c_{k}\}$. If $F(F(t(x),t(y)),t(z))<t(1)$, then $F(F(t(x),t(y)),t(z))=F(t(x), F(t(y),t(z)))\in M\setminus C$ since $F(t(y),t(z))\in M\setminus C$.
This implies $F(F(t(x),t(y)),t(z))=F(b_{k}, t(z))$. Since $t(z) < \infty$ and $F$ is strict, we have $F(t(x), t(y))=b_{k} $, contrary to $F(t(x),t(y))\in [b_{k}, d_{k}]\setminus\{c_{k}\}$. Therefore, $t(1) \leq \min\{F(b_{k}, t(z)), F(F(t(x),t(y)),t(z))\}<\infty$, i.e.,
$F((b_{k}, F(t(x), t(y))), t(z)) \subseteq (t(1),\infty]$.

Thus by the choice of $t(z)$, we have
$$F((b_{k}, F(t(x), t(y))),(M\setminus C)\setminus\{\infty\}) \subseteq (t(1),\infty]$$
since $F$ is strict. This follows from the definition of $H_{k}$ that
$$F(\cup_{k\in K_{1}}H_{k},(M\setminus C)\setminus\{\infty\})\subseteq [t(1),\infty],$$
which results in $F(\cup_{k\in K_{1}}H_{k},M\setminus\{\infty\})\cap (M\setminus C)=\emptyset$ since $t\mid_{[0, \alpha]}$ is strictly increasing. Therefore, $F(\cup_{k\in K_{1}}H_{k}, M)\cap (M\setminus t(0))= \emptyset$ since $t(1)<\infty$.

(ii)$\Rightarrow$(i) Suppose that $F(\cup_{k\in K_{1}}H_{k}, M)\cap (M\setminus t(0))= \emptyset$. It is enough to prove that $T(T(x,y),z)=T(x, T(y,x))$ for arbitrary $x,y,z\in [0,1]$. We first prove that $T(T(x,y),z)=t^{[-1]}F((F(t(x),t(y))),t(z))$ for arbitrary $x,y,z\in [0,1]$. In fact, from Eq.~\eqref{eq:12} we have $T(x,y)=t^{[-1]}(F(t(x),t(y)))$. We consider two cases as below.

(1) If $F(t(x),t(y))\in M$, then $t\circ t^{[-1]}(F(t(x),t(y)))=F(t(x),t(y))$ implies
\begin{align}
 T(T(x,y),z) &= t^{[-1]}F(t\circ t^{[-1]}(F(t(x),t(y))),t(z)) \nonumber \\
 &=  t^{[-1]}F((F(t(x),t(y))),t(z)).\nonumber
\end{align}

(2) If $F(t(x),t(y))\notin M$, then $F(t(x),t(y))\in [0,t(0)]$ or there exists a $k \in K$ such that $F(t(x),t(y))\in [b_{k}, d_{k}]\setminus\{c_{k}\}$. If $F(t(x),t(y))\in [0,t(0)]$ then $(F(t(x),t(y)),t(0)]\in H_{\kappa}$, thus $$F((F(t(x),t(y)), t(0)], t(z))\cap (M\setminus t(0))= \emptyset.$$
It follows from Lemma \ref{lem3.5} that
$$ t^{[-1]}(F(F(t(x),t(y))),t(z))= t^{[-1]}(F(t(0),t(z))).$$
Therefore,
$$T(T(x, y),z)=T(0,z)=t^{[-1]}F(t(0),t(z))=t^{[-1]}F((F(t(x),t(y))),t(z))$$
since $T(x, y)= t^{[-1]}(F(t(x),t(y)))=0$.

If there exists a $k \in K$ such that $F(t(x),t(y))\in [b_{k}, d_{k}]\setminus{c_{k}}$ then $(b_{k}, T(t(x),t(y))] \in H_{k}$, thus
$$F((b_{k}, T(t(x),t(y))] , t(z))\cap (M\setminus t(0))= \emptyset.$$
It follows from Lemma \ref{lem3.5} that
$$ t^{[-1]}(F(F(t(x),t(y))),t(z))=  t^{[-1]}(F(b_{k},t(z))).$$
Therefore,
$$T(T(x,y),z)=t^{[-1]}F(t\circ t^{[-1]}F(t(x),t(y)),t(z))=t^{[-1]}F(b_{k},t(z))=t^{[-1]}F((F(t(x),t(y))),t(z))$$
since $b_{k}=t\circ t^{[-1]}(F(t(x),t(y)))$.

Analogously, we have $T(x, T(y,z))=t^{[-1]}F(t(x),F(t(y),t(z)))$ for all $x,y,z\in [0,1]$.

Therefore, $T(T(x,y),z)=T(x, T(y,x))$ since $F$ is associative.
\end{proof}

\subsection{Continuity}

This subsection is devoted to analyzing the relationship between the continuity of both $t$ and $T$ given by Eq. \eqref{eq:12}.

Let $t:[0,1]\rightarrow [0,\infty]$ be a left continuous non-decreasing function, $F\in \varGamma$ and the function $T$ be given by Eq. \eqref{eq:12}. We write
$$ T(x_{0}^{-},y_{0}^{-})=\lim_ {x \nearrow x_{0}, y \nearrow y_{0}} T(x, y),$$
 $$ T(x_{0}^{+},y_{0}^{+})=\lim_ {x \searrow x_{0}, y \searrow y_{0}} T(x, y).$$
Then $T$ is left continuous (resp. right continuous) at $(x_0,y_0)$ if and only if  $T(x_{0}^{-},y_{0}^{-})=T(x_0,y_0)$ (resp. $T(x_{0}^{+},y_{0}^{+})=T(x_0,y_0)$); $T$ is continuous at $(x_0,y_0)$ if and only if it is left continuous and right continuous at $(x_0,y_0)$.

In what follows, we always suppose that
 $F\in \varGamma$ is continuous. Then the following proposition is obvious.
 \begin{proposition}\label{prop3.15}
Let $F\in \varGamma$ be continuous. If $t:[0,1]\rightarrow [0,\infty]$ is a continuous non-decreasing function, then the function $T$ given by Eq.~\eqref{eq:12} is right continuous.
\end{proposition}
%\begin{proof}If $t:[0,1]\rightarrow [0,\infty]$ is a continuous non-decreasing function, then $t^{[-1]}$ is a right continuous non-decreasing function. This implies that $T$ is right continuous since $F$ is continuous.
%\end{proof}
\begin{proposition}\label{prop3.16}
Let the function $T$ be given by Eq.~\eqref{eq:12}$, F\in \varGamma$ be continuous and $t:[0,1]\rightarrow [0,\infty]$ be a left continuous non-decreasing function. For any $(x,y)\in[0,1]^{2}$, if  $F(t(x),t(y))\in [0,\infty]\setminus\mathbb{H}$, then $T$ is left continuous at $(x,y)$. Furthermore, if $t$ is continuous at $(x,y)$, then $T$ is continuous at $(x,y)$.
\end{proposition}
\begin{proof} For any $(x,y)\in[0,1]^{2}$, if $F(t(x),t(y))\in [0,\infty]\setminus\mathbb{H}$, then there exists an $\varepsilon >0$ such that $t^{[-1]}$ is continuous on
$[F(t(x),t(y))-\varepsilon, F(t(x),t(y))+\varepsilon]$. Thus $T(x^{-},y^{-})=t^{[-1]}(F(t(x^-),t(y^-)))=t^{[-1]}(F(t(x),t(y)))=T(x,y)$ since $t$ is left continuous and $F$ is continuous. Moreover, if $t$ is continuous at $(x,y)$, then by Proposition \ref{prop3.15} we have $T(x^{-},y^{-})=T(x,y)=T(x^{+},y^{+})$.
\end{proof}

Trivially, if $t:[0,1]\rightarrow [0,\infty]$ is strictly increasing, then $F(t(x),t(y))\in [0,\infty]\setminus\mathbb{H}$ for all $(x,y)\in[0,1]^{2}$. Therefore, from Proposition \ref{prop3.16}, we have the following corollary.
\begin{corollary}\label{cor3.3}
Let the function $T$ be given by Eq.~\eqref{eq:12} and $F\in \varGamma$ be continuous.
\renewcommand{\labelenumi}{(\roman{enumi})}
\begin{enumerate}
\item If $t:[0,1]\rightarrow [0,\infty]$ is a left continuous strictly increasing function, then $T$ is left continuous.
\item If $t:[0,1]\rightarrow [0,\infty]$ is a continuous strictly increasing function, then $T$ is continuous.
\end{enumerate}
\end{corollary}

Generally, the converse implications of both Corollary \ref{cor3.3} (i) and (ii) do not hold.
\begin{example}\label{exp5.1}
\emph{Let the function $t:[0,1]\rightarrow [0,\infty]$ be given by
\begin{equation*}
 		t(x)=\begin{cases}
 			4 & \ \hbox{if } x\in[0,0.5],\\
 			2x+4 &  \ \hbox{if } x\in(0.5,1]\\
 		\end{cases}
 	\end{equation*}
 and the function $F:[0,\infty]^{2}\to[0,\infty]$ be defined by $F(x,y)=x+ y$. Then by Eq.~\eqref{eq:12},
 $$T(x, y) = 1 \hbox{ ~for any } (x,y)\in[0,1]^2.$$
It is easy to see that $T$ is continuous on $[0,1]^2$. However,  $t:[0,1]\rightarrow [0,\infty]$ is neither strictly increasing nor continuous.}
\end{example}

Fortunately, we have the following proposition.
\begin{proposition}\label{prop3.17}
Let $F\in \varGamma$ be continuous and $t:[0,1]\rightarrow [0,\infty]$ be a left continuous strictly increasing function. Then
 the function $T$  given by Eq.~\eqref{eq:12} is continuous at $(x,y)$ if and only if  $M\cap[F(t(x),t(y)),F(t(x^+),t(y^+))]$ is at most a one-element set.
\end{proposition}
\begin{proof}Let $t:[0,1]\rightarrow [0,\infty]$ be a strictly increasing function. Then, from (i) of Corollary \ref{cor3.3}, $T$ is left continuous. Therefore, $T$ is right continuous at $(x,y)\in[0,1]^2$ if and only if $t^{[-1]}(F(t(x),t(y)))=t^{(-1)}(F(t(x^+),t(y^+)))$ if and only if $M\cap[F(t(x),t(y)),F(t(x^+),t(y^+))]$ is at most a one-element set.
\end{proof}

Note that $M\cap[F(t(x),t(1)),F(t(x^+),t(1^+))]$ is always at most a one-element set. Therefore, from Proposition \ref{prop3.17} we have the following corollary.
\begin{corollary}\label{coro5.1}
Let $F\in \varGamma$ be continuous and $t:[0,1]\rightarrow [0,\infty]$ be a left continuous strictly increasing function. The following are equivalent:
\renewcommand{\labelenumi}{(\roman{enumi})}
\begin{enumerate}
\item The function $T$ given by Eq.~\eqref{eq:12} is continuous.
\item $M\cap[F(t(x),t(y)),F(t(x^+),t(y^+))]$ is at most a one-element set for all $x,y\in[0,1)$.
\end{enumerate}
\end{corollary}

 \begin{proposition}\label{prop3.18}
 Let $ F\in \varGamma$ be continuous and $t:[0,1]\rightarrow [0,\infty]$ be a left continuous non-decreasing function. If there exists an $x_0\in[0,1)$ such that $t^{[-1]}(t(x_{0}))=t^{[-1]}(t(x_{0}^+))$ and $F(t(x_{0}^+),t(x_{0}^+))=t(x_{0}^+)$, then the function $T$ given by Eq.~\eqref{eq:12} is right continuous at $(x_0,y)$ for all $y\leq x_0$.
 \end{proposition}
 \begin{proof}Let $y\leq x_0$. Then $t(y)\leq t(y^{+})\leq t(x_{0})\leq t(x_{0}^{+})$  by the monotonicity of $t$. Hence, by the monotonicity of $F$, $$t(x_{0}^+)=F(t(x_{0}^+),t(x_{0}^+))\geq F(t(x_{0}),t(y))\geq t(x_{0}).$$ Meanwhile, by Proposition \ref{pr1.6}, $F(t(x_{0}^+),t(y^+))=t(x_{0}^+)$ since $F(t(x_{0}^+),t(x_{0}^+))=t(x_{0}^+)$. Thus $$[F(t(x_{0}),t(y)),F(t(x_{0}^+),t(y^+))]\subseteq[t(x_{0}),t(x_{0}^{+})].$$ Moreover, $t^{[-1]}(F(t(x_0),t(y)))=t^{[-1]}(F(t(x_0^+),t(y^+)))$ since $t^{[-1]}(t(x_{0}))=t^{[-1]}(t(x_{0}^+))$, i.e., $T(x_0,y)=T(x_0^{+},y^{+})$. Therefore, $T$ is right continuous at $(x_0,y)$.
 \end{proof}

Notice that the condition $t^{[-1]}(t(x_{0}))=t^{[-1]}(t(x_{0}^+))$ in Proposition \ref{prop3.18} cannot be deleted generally.
\begin{example}\label{ex5.1.05}
\emph{Let the function $t:[0,1]\rightarrow [0,\infty]$ be given by
\begin{equation*}
 		t(x)=\begin{cases}
 			5x& \ \hbox{if } x\in[0,0.4],\\
 			2& \ \hbox{if } x\in(0.4,0.5],\\
 			4 & \ \hbox{if } x\in(0.5,1]\\
 		\end{cases}
 	\end{equation*}
and the function $F:[0,\infty]^{2}\to[0,\infty]$ be defined by $F(x,y)=\max\{x,y\}$. Then by Eq.~\eqref{eq:12},
\begin{equation*}
 		T(x,y)=\begin{cases}
 			\max\{x,y\} & \ \hbox{if } x+y\leq0.4,\\
 			0.5 &  \ \hbox{if } (x,y)\in [0,0.5]^2 \setminus\{x+y\leq 0.4\},\\
 			1 & \ \mbox{otherwise}.
 		\end{cases}
 	\end{equation*}
It is easy to check that  $F(t(x_{0}^+),t(x_{0}^+))=t(x_{0}^+)$ for all $x_0\in[0,1)$, $t^{[-1]}(t(0.5))=0.5$ and  $t^{[-1]}(t(0.5^+))=1$. However, $T(0.5,0.5)=0.5$ and $T(0.5^+,0.5^+)=1$, i.e.,
$T$ is not right continuous at $(0.5,0.5)$.}
\end{example}

\begin{corollary}\label{cor7.8}
Let $ F\in \varGamma$ be continuous and $t:[0,1]\rightarrow [0,\infty]$ be a left continuous strictly increasing function. If there exists an $x_0\in[0,1)$ such that $F(t(x_{0}^+),t(x_{0}^+))=t(x_{0}^+)$, then the function $T$ given by Eq.~\eqref{eq:12} is continuous at $(x_0,y)$ for all $y\in[0,1]$.
\end{corollary}
\begin{proof}Let $ F\in \varGamma$ be continuous and $t:[0,1]\rightarrow [0,\infty]$ be a left continuous strictly increasing function. Then, from (i) of Corollary \ref{cor3.3}, $T$ is left continuous. Further, by Proposition \ref{prop3.18}, $T$ is continuous at $(x_0,y)$ for all $y\leq x_0$. On the other hand, for all $y>x_0$, we have $t(x_{0})\leq t(x_{0}^{+})\leq t(y)\leq t(y^{+})$ by the monotonicity of $t$. Because of $ F\in \varGamma$ being continuous and $F(t(x_{0}^+),t(x_{0}^+))=t(x_{0}^+)$, from Proposition \ref{pr1.6} $F(t(x_{0}^+),t(y)) = t(y)\mbox{ and } F(t(x_{0}^+),t(y^+)) = t(y^+)$. These imply that $$t(y)\leq F(t(x_{0}),t(y))\leq F(t(x_{0}^+),t(y))=t(y),$$ i.e., $F(t(x_{0}),t(y))=t(y)$. Thus $$[F(t(x_{0}),t(y)),F(t(x_{0}^+),t(y^{+}))]=[t(y),t(y^{+})],$$ and we have $t^{[-1]}(F(t(x_0),t(y)))=t^{[-1]}(F(t(x_0^+),t(y^+)))$, i.e., $T(x_0,y)=T(x_0^+,y^+)$. Therefore, $T$ is continuous at $(x_0,y)$ for all $y\in[0,1]$.
  \end{proof}

\begin{lemma}\label{le7.3}
Let $ F\in \varGamma$ be continuous. If $t:[0,1]\rightarrow [0,\infty]$ is a left continuous strictly increasing function and $F(t(0^+),t(0^+))=t(0^+)$, then the function $T$ given by Eq.~\eqref{eq:12} is continuous.
\end{lemma}
\begin{proof}From Corollary \ref{cor7.8}, we have that $T$ is continuous on $(0,1]^2$. Below we show that $T$  is continuous at $(0,x)$ for all $x\in[0,1]$. Because $ F\in \varGamma$ is continuous and $F(t(0^+),t(0^+))=t(0^+)$, by Proposition \ref{pr1.6}, $T$ is continuous at $(0,x)$ for all $x\in(0,1]$. On the other hand, $t(0)\leq F(t(0),t(0))\leq F(t(0^+),t(0^+))=t(0^+)$, which means
 $F(t(0),t(0))=F(t(0^+),t(0^+))$. Thus we have $T(0,0)=T(0^+,0^+)$. Therefore, $T$ is continuous.
 \end{proof}

%Note that, from Example \ref{exp5.1}, the converse of Lemma \ref{le7.3} does not generally hold.

\begin{lemma}\label{le7.4}
Let $ F\in \varGamma$ be continuous and $t:[0,1]\rightarrow [0,\infty]$ be a left continuous non-decreasing function. If the function $T$ given by Eq.~\eqref{eq:12} is continuous at $(0,0)$ and $T(0,0)=0$, then $t(0)\in M\setminus \mathbb{D}$ and $F(t(0^+),t(0^+))=t(0^+)$.
\end{lemma}
\begin{proof}Let the function $T$ given by Eq.~\eqref{eq:12} be continuous at $(0,0)$ and $T(0,0)=0$. Then $$0=T(0,0)=T(0^{+},0^{+})=t^{[-1]}(F(t(0^+),t(0^+))).$$ If $t(0)\notin M\setminus \mathbb{D}$, then $t(0)\in \mathbb{D}$. In this case, we have $t^{[-1]}(t(0))>0$, thus, $$T(0^{+},0^{+})=t^{[-1]}(F(t(0^+),t(0^+)))\geq t^{[-1]}t(0^+))\geq t^{[-1]}(t(0))>0,$$ a contradiction.
 	Thus $t(0)\in M\setminus \mathbb{D}$, which implies $t^{[-1]}(t(0))=0$. Then $t^{[-1]}(F(t(0^+),t(0^+)))=t^{[-1]}(t(0))$, hence $$M\cap[t(0), F(t(0^+),t(0^+))]=\{t(0)\}.$$  The last equality implies that $F(t(0^+), t(0^+))=t(0^+)$ since  $F(t(0^+),t(0^+))\geq t(0^+)$.
 \end{proof}

 We finally have the following theorem.
\begin{theorem}\label{Thm3.19}
Let $ F\in \varGamma$ be continuous and $t:[0,1]\rightarrow [0,\infty]$ be a left continuous non-decreasing function. If $F(t(x),t(0^+))<F(t(x),t(y))$ for all $x,y\in(0,1)$, then the following are  equivalent:
\renewcommand{\labelenumi}{(\roman{enumi})}
\begin{enumerate}
\item The function $T$ given by Eq.~\eqref{eq:12} is a continuous t-conorm.
\item $t$ is a continuous strictly increasing function on $(0,1]$ and $F(t(0^+),t(0^+))=t(0^+)$.
\end{enumerate}
\end{theorem}
\begin{proof}
	(i)$\Rightarrow$(ii) Let the function $T$ given by Eq.~\eqref{eq:12} be a continuous t-conorm. If $t$ is not strictly increasing, then there exist two elements $x,y\in[0,1]$ with $x<y$ such that $t(x)=t(y)$ and $t^{[-1]}(t(x))=t^{[-1]}(t(y))=y$. Thus we have $$T(x,0)=t^{[-1]}(F(t(x),t(0)))\geq t^{[-1]}(F(t(x),0))\geq t^{[-1]}(t(x))=y>x,$$ a contradiction. Therefore, $t$ is strictly increasing. Meanwhile, from Lemma \ref{le7.4}, we have $F(t(0^+),t(0^+))=t(0^+)$. Thus, by Proposition \ref{pr1.6},  $$F(t(x),t(0^+))=t(x)\leq t(x^{+}) \mbox{ for all } x\in(0,1).$$ If there exists an $x\in(0,1)$ such that  $t(x)<t(x^{+})$, then there is a $z\in(0,1)$ near enough to 0 such that $F(t(x),t(z))<t(x^+)$ by the continuity of $F$. Thus
 	 $$F(t(x),t(z))<t(x^+)=F(t(x^+),t(0^+))<F(t(x^{+}),t(z))\leq F(t(x^{+}),t(z^{+})),$$
 	which yields that $t^{[-1]}(F(t(x),t(y)))\neq t^{[-1]}(F(t(x^{+}),t(y^{+})))$, i.e., $T(x,y)\neq T(x^{+},y^{+})$, contrary to the fact that $T$ is continuous. Therefore, $t(x)= t(x^{+}) \mbox{ for all } x\in(0,1),$ i.e., $t$ is continuous.
 	
 	(ii)$\Rightarrow$(i) From Lemma \ref{le7.3},  $T$ is continuous. Below we show that $T$ is a t-conorms. Let $t$ be a continuous strictly increasing function on $(0,1]$. Then, obviously, $T$ satisfies the monotonicity, the associativity and the commutativity. Because $F\in\varGamma$ is continuous and $F(t(0^+),t(0^+))=t(0^+)$, by Proposition \ref{pr1.6}, we have $F(t(x),t(0^+))=t(x)$ for all $x\in(0,1]$. Thus
 	$$t(x)\leq F(t(x),t(0))\leq F(t(x),t(0^+))=t(x) \mbox{ for all } x\in(0,1],$$
 	i.e., $F(t(x),t(0))=t(x)$ for all $x\in(0,1]$.  On the other hand, $t(0^+)\leq F(t(0),t(0^+))\leq F(t(0^+),t(0^+))=t(0^+),$
 	i.e., $F(t(0),t(0^+))=t(0^+)$, thus $F(t(x),t(0))=t(x)$ for all $x\in[0,1]$.
 	Then by Eq.~\eqref{eq:12} we get $$T(x,0)=t^{[-1]}(F(t(x),t(0)))=t^{[-1]}(t(x))=x$$ for any $x\in[0,1]$. Therefore, by Definition \ref{de2.2} $T$ is a t-conorm.
\end{proof}

Notice that the condition $F(t(x),t(0^+))<F(t(x),t(y))$ for all $x,y\in(0,1)$ in Theorem \ref{Thm3.19} cannot be dropped generally.
\begin{example}\label{ex5.1.5}
 	\emph{Let the function $t:[0,1]\to[0,\infty]$ be given by\begin{equation*}
         t(x)=\begin{cases}
              2x& \hbox{if }\ x\in[0,\frac{1}{2}],\\
              4x& \hbox{if }\ x\in(\frac{1}{2},1]\\
              \end{cases}
         \end{equation*}
 	and the function $F:[0,\infty]^{2}\to[0,\infty]$ be given by $F(x,y)=\max\{x,y\}$.
 	Then by Eq~\eqref{eq:12},
 \begin{equation*}
 		T(x,y)=\max\{x,y\} \mbox{ for all } (x,y)\in[0,1]^2.
 	\end{equation*}
 	It is easy to check that $T$ is a continuous t-conorm. However, $t$ is not a continuous function on $(0,1]$.}
 \end{example}

% \begin{proposition}\label{pr5.1.20}
% 	Let a function $T$ be given by Eq.~\eqref{eq:12}$,  F\in \varGamma$ be strict and $t:[0,1]\rightarrow [0,\infty]$ be a left continuous non-decreasing function. If $T(0,0)=0$, then the following are  equivalent:
% \renewcommand{\labelenumi}{(\roman{enumi})}
% \begin{enumerate}
% \item $T$  is  continuous.
%\item $t:[0,1]\rightarrow[0,\infty]$ is a continuous strictly increasing  function.
% \end{enumerate}
% \end{proposition}
% \begin{proof} From  Corollary \ref{cor3.3}, (ii)$\Rightarrow$(i) is obvious. Below we show that (i)$\Rightarrow$(ii).
% 	Becaues of $ F\in \varGamma$ is strict, then $F(t(x),t(0^+))<F(t(x),t(y))$ for any $x,y\in(0,1)$.
% 	By Proposition \ref{prop3.19}, $t$ is a strictly increasing continuous function on $(0,1]$ and  $F(t(0^+),t(0^+))=t(0^+)$. We get
% 	$$t(0^+)\leq F(t(0),t(0^+))\leq F(t(0^+),t(0^+))=t(0^+),$$
% 	i.e., $F(t(0),t(0^+))=t(0^+)$, then $t(0)=t(0^+)$ since $ F$ is strict. Thus, $t:[0,1]\rightarrow[0,\infty]$ is a strictly increasing continuous function.
% \end{proof}

\begin{corollary}\label{cor7.9}
 	Let $ F\in \varGamma$ be strict and $t:[0,1]\rightarrow [0,\infty]$ be a left continuous non-decreasing function. Then the following are equivalent:
 \renewcommand{\labelenumi}{(\roman{enumi})}
 \begin{enumerate}
 \item The function $T$ given by Eq.~\eqref{eq:12} is a continuous t-conorm.
\item $t:[0,1]\rightarrow[0,\infty]$ is a continuous strictly increasing function and $F(t(0),t(0))=t(0)$.
 \end{enumerate}
 \end{corollary}
\begin{proof} From  Theorem \ref{Thm3.19}, (ii)$\Rightarrow$(i) is obvious. Below we show that (i)$\Rightarrow$(ii).
 	Becaues of $ F\in \varGamma$ being strict, $F(t(x),t(0^+))<F(t(x),t(y))$ for any $x,y\in(0,1)$.
 	Then by Theorem \ref{Thm3.19}, $t$ is a strictly increasing continuous function on $(0,1]$ and  $F(t(0^+),t(0^+))=t(0^+)$. Thus
 	$$t(0^+)\leq F(t(0),t(0^+))\leq F(t(0^+),t(0^+))=t(0^+),$$
 	i.e., $F(t(0),t(0^+))=t(0^+)$. Thus $t(0)=t(0^+)$ since $ F$ is strict. Therefore, $t:[0,1]\rightarrow[0,\infty]$ is a strictly increasing continuous function and $F(t(0),t(0))=t(0)$.
 \end{proof}
% \begin{example}\label{ex5.1.5}
% 	\emph{Let the function $t:[0,1]\to[0,\infty]$ be given by  $$t(x)=2x \mbox{ ~for all } x\in[0,1]$$
% 	 and function $F:[0,\infty]^{2}\to[0,\infty]$ be given by
% 	\begin{equation*}
% 		F(x,y)= x+y+xy \mbox{ ~for all } (x,y)\in[0,\infty]^2.
% 	\end{equation*}
% 	Then by Eq~\eqref{eq:12}, \begin{equation*}
% 		T(x,y)=\begin{cases}
% 			x+y+2xy & \hbox{if } \ x+y+2xy\leq1,\\
% 			1 &\ \mbox{otherwise}.
% 		\end{cases}
% 	\end{equation*}
% 	It is easy to check that $ F\in \varGamma$ is strict which satisfies $F(t(0),t(0))=t(0)$. From Corollary \ref{cor7.9}, $T$ is a continuous t-conorm.}
% \end{example}

\section{Conclusions}

One can easily check that our results are suitable for all left continuous non-increasing functions also. So the main contributions of this article include that we gave the concept of a weak pseudo-inverse of a monotone function for overcoming the difficulty that a function $T:[0,1]^2\rightarrow[0,1]$ defined by Eq.(\ref{eq:2}) isn't associative when $t$ is a left continuous monotone function, answered what is the characterization of a left continuous monotone function $t: [0,1]\rightarrow [0,\infty]$ such that the function given by Eq.(\ref{eq:12}) is associative, and furthermore presented the idempotence, the limit property, the conditional cancellation law and the continuity of the function $T$  given by Eq.(\ref{eq:12}), respectively. It is regrettable that generally, our results aren't true when $t$ is a right continuous monotone function. For instance, let $F(x,y)=\max\{x,y\}$ and the function $t:[0,1]\rightarrow [0,1]$ be defined by \begin{equation*}
 t(x)=\begin{cases}
x & \hbox{if }\ x\in[0,\frac{1}{2}],\\
\frac{1}{2} &  \hbox{if }\ x\in(\frac{1}{2},\frac{3}{4}),\\
\frac{3}{4} & \hbox{if }\ x\in[\frac{3}{4},1].
\end{cases}
\end{equation*}
It is easy to see that $t$ is a right continuous non-decreasing function. Then $M=[0,\frac{1}{2}]\cup\{\frac{3}{4}\}$. Thus from Definition 5.2 of \cite{YM2024}, we have $I(M)=\emptyset$, so that $I(M) \cap (M\setminus \{t(1)\})=\emptyset$. Therefore, by Corollary 5.3 of \cite{YM2024} we know that the following function $T:[0,1]^2\rightarrow[0,1]$ given by Eq.(6) of \cite{YM2024} is associative since $t(1)<1$:
\begin{equation*}
T(x,y)=\left\{
  \begin{array}{ll}
  \max\{x,y\} & \hbox{if }(x,y)\in[0,\frac{1}{2})^{2}, \\
   \frac{3}{4}  & \hbox{otherwise.}
  \end{array}
\right.
\end{equation*}

On the other hand, from Definition \ref{def6.1} we have $\mathfrak{T}(M)=\emptyset$. This follows that $\mathfrak{T}(M) \cap (M\setminus \{t(0)\})=\emptyset$. However, from Eq.(\ref{eq:12}),
\begin{equation*}
T(x,y)=\left\{
  \begin{array}{ll}
  \max\{x,y\} & \hbox{if }(x,y)\in[0,\frac{1}{2})^{2}, \\
   \frac{3}{4} & \hbox{if }(x,y)\in[0,\frac{3}{4})^{2}\setminus [0,\frac{1}{2})^{2}, \\
   1  & \hbox{otherwise.}
  \end{array}
\right.
\end{equation*}
Put $x=\frac{1}{4}, y=\frac{1}{4}, z=\frac{1}{2}$. Then $T(T(x,y),z)=T(T(\frac{1}{4},\frac{1}{4}),\frac{1}{2})=T(\frac{1}{4},\frac{1}{2})=\frac{3}{4}$, $T(x,T(y,z))=T(\frac{1}{4},T(\frac{1}{4},\frac{1}{2}))=T(\frac{1}{4},\frac{3}{4})=1$. Thus, $T$ isn't an associative function.

%\begin{figure}[htbp]
%	\centering
%\subfigure[] {
%\includegraphics[width=6cm]{T2.eps}
%	}
%	\quad
%	\label{T1}
%\subfigure[]{
%		\includegraphics[width=6cm]{T1.eps}
%	}
%	\label{T2}
%\end{figure}
%
%
%\begin{figure}[htbp]
%	\centering
%	\subfigure[]{
%		\includegraphics[width=6cm]{T3.eps}
%	}
%	\quad
%	\label{T2}
%\subfigure[]{
%		\includegraphics[width=6cm]{T4.eps}
%	}
%\caption{(a) a strictly increasing function; (b) a right continuous non-decreasing function; (c) a kind of non-decreasing functions; (d) a left continuous non-decreasing function.}
%\end{figure}

\end{document}